\documentclass[12pt]{amsart}

\newcommand {\supplus}{\mathop{{\supset}\llap{\raise 
0.5pt\hbox{\normalfont\small+}\hskip 0.5pt}}} 

\newcommand {\subplus}{\mathop{{\subset}\llap{\raise 
0.5pt\hbox{\normalfont\small+}\hskip 0.5pt}}}  

\newcommand {\Cee}    {{\mathbb  C}}

\newcommand {\Kee}    {{\mathbb  K}}
\newcommand {\Nee}    {{\mathbb  N}}

\newcommand {\Qee}    {{\mathbb  Q}}
\newcommand {\Ree}    {{\mathbb  R}}

\newcommand {\Zee}    {{\mathbb  Z}}

\newcommand {\fa}     {{\mathfrak{a}}}

\newcommand {\faut}   {{\mathfrak{aut}}}

\newcommand {\fder}   {{\mathfrak{der}}}   %
\newcommand {\fdiff}  {{\mathfrak{diff}}}

\newcommand {\fe}     {{\mathfrak{e}}}

\newcommand {\ff}     {{\mathfrak{f}}}
\newcommand {\fF}     {{\mathfrak{F}}}
\newcommand {\fg}     {{\mathfrak{g}}}    %
\newcommand {\fG}     {{\mathfrak{G}}}    %
\newcommand {\fgl}    {{\mathfrak{gl}}}  %
\newcommand {\fh}     {{\mathfrak{h}}}

\newcommand {\fI}    {{\mathfrak{I}}}    %
\newcommand {\fk}     {{\mathfrak{k}}}

\newcommand {\fn}     {{\mathfrak{n}}}
\newcommand {\fns}     {{\mathfrak{ns}}}
\newcommand {\fo}     {{\mathfrak{o}}}
\newcommand {\fosp}   {{\mathfrak{osp}}}
\newcommand {\fp}    {{\mathfrak{p}}}   %

\newcommand {\fpo}    {{\mathfrak{po}}}

\newcommand {\fr}     {{\mathfrak{r}}}

\newcommand {\fs}     {{\mathfrak{s}}}

\newcommand {\fsl}    {{\mathfrak{sl}}}

\newcommand {\fsp}    {{\mathfrak{sp}}}

\newcommand {\fsq}    {{\mathfrak{sq}}}

\newcommand {\fvect}  {{\mathfrak{vect}}}   %

\newcommand {\cal} {\mathcal}

\newcommand {\cF}     {{\cal F}}

\newcommand {\cL}     {{\cal L}}
\newcommand {\cM}     {{\cal M}}

\newcommand {\cO}     {{\cal O}}

\newcommand {\cX}     {{\cal X}}

\def \opname#1#2%
  {\expandafter\newcommand \csname #1\endcsname {{\mathop{#2}\nolimits}}}

\newcommand{\rmname}[1]
  {\expandafter\newcommand \csname #1\endcsname {{\operatorname{#1}}}}

\newcommand{\rmnameii}[2]
  {\expandafter\newcommand \csname #1\endcsname {{\operatorname{#2}}}}

\rmname{act}
\rmname{Ad}
\rmname{Add}
\rmname{ad}
\rmname{Alt}
\rmname{alt}
\rmname{Ann}
\rmname{antidiag}
\rmname{Ber}
\rmname{ber}
\rmname{Br}
\rmname{card}
\rmname{ch}
\rmname{Char}
\rmname{cem}
\rmname{cj}
\rmname{Cliff}
\rmname{cntr}
\rmname{codim}
\rmname{coind}
\rmname{const}
\rmname{col}
\rmname{cork}
\rmname{cpr}
\rmname{diag}
\rmnameii{Div}{div}
\rmname{Def}
\rmname{Der}
\rmname{Diff}
\rmname{Dim}
\rmname{End}
\rmname{Even}
\rmname{Ext}
\rmname{gr}
\rmname{Hom}
\rmname{HT}
\rmnameii{Ht}{ht}
\rmname{hwt}
\rmname{Id}
\rmname{id}
\rmname{ind}
\rmname{Ind}
\rmname{Inf}
\rmname{irr}
\rmname{Le}
\rmname{Lie}
\rmname{lwt}
\rmname{mult}
\rmname{Mat}
\rmname{Mor}
\rmname{nm}
\rmname{Ob}
\rmname{Odd}
\rmname{Osc}
\rmname{per}
\rmname{Pic}
\rmname{pr}
\rmname{pro}
\rmname{Prime}
\rmname{Proj}
\rmname{prt}
\rmname{pt}
\rmname{Q}
\rmname{qet}
\rmname{qtr}
\rmname{rd}
\rmname{rk}
\rmname{row}
\rmname{Res}
\rmname{salt}
\rmname{Sch}
\rmname{SBr}
\rmname{scalar}
\rmname{Ser}
\rmname{sign}
\rmname{Smbl}
\rmname{spin}
\rmname{ssym}
\rmname{str}
\rmname{st}
\rmname{sgn}
\rmname{sq}
\rmname{symm}
\rmname{supp}
\rmname{Supp}
\rmname{St}
\rmname{Spec}
\rmname{Spm}
\rmname{tr}
\rmname{vpt}
\rmname{weyl}
\rmname{Weyl}
\rmname{Witt}

\opname{vvol}  {{v\hspace{-0.1ex}o\hspace{-0.02ex}l\/}}
\opname{pnt}  {\text{\normalfont pt}}
\opname{Span} {{Span}}
\opname{slim} {\overline{\lim}}
\opname{Vol}  {{V\hspace{-0.55ex}o\hspace{-0.02ex}l\/}}
\opname{Par}  {{P\hspace{-0.3ex}a\hspace{-0.05ex}r\/}}

\newcommand {\ev} {{\bar0}}
\newcommand {\od} {{\bar1}}

\newcommand {\tto} {\longrightarrow}
\newcommand {\pder}[1] {{\frac{\partial}{\partial {#1}}}}

\newcommand {\bcdot}   {\mathbin{\hbox{\raise.4ex\hbox{\bf.}}}} 

\newcommand {\secno} {}
\newcommand {\ssecfont} {\normalfont\bf}

\newtheorem{Theorem}{\secno Theorem}

\newtheorem{Proposition}[Theorem]{\secno Proposition}

\newtheorem{Problem}[Theorem]{\secno Problem}

\newenvironment {th*}[1]
    {\gdef\thname{#1} \begin{thn}}%
    {\end{thn}}
\newtheorem{thn}[Theorem] {\thname}

\theoremstyle{definition}

\newenvironment {ex*}[1]
    {\gdef\thname{#1} \begin{exn}}%
    {\end{exn}}
\newtheorem{exn}[Theorem]{\thname}

\theoremstyle{remark}

\newtheorem{Remark}[Theorem]{\secno Remark}

\newenvironment {rem*}[1]
    {\gdef\thname{#1} \begin{remn}}%
    {\end{remn}}
\newtheorem{remn}[Theorem]{\thname}

\newcommand {\ssec}{\subsection*}

\newcommand {\ssbegin}[2]
  {\def \secno {\gdef \secno {}{\ssecfont #1. }}%
   \begin{#2}}
\begin{document}

\title[Lie superalgebras of supermatrices of complex 
size]{Lie superalgebras of supermatrices of complex 
size. Their generalizations and related integrable systems}

\author{Pavel Grozman, Dimitry Leites} 

\address{Dept.  of Math., Univ.  of Stockholm, Roslagsv.  101, 
Kr\"aftriket hus 6, S-106 91, Stockholm, 
Sweden\\e-mail:mleites\@matematik.su.se}

\thanks{We are thankful to V.~Kornyak (JINR, Dubna) who checked the 
generators and relations with an independent program and compared 
convenience of the Serre relations with that of our ones.   Financial
support of the Swedish Institute and NFR is gratefully  acknowledged. 
We are thankful to B.~Feigin and Shi Kangjie-laoshi for  their shrewd
questions and to S.~Shnider, G.~Post and M.~Vasiliev for the timely 
information.}



\begin{abstract} We distinguish a class of simple filtered Lie 
algebras $LU_\fg(\lambda)$ of polynomial growth with increasing 
filtration and whose associated graded Lie algebras are not simple.  
We describe presentations of such algebras.  The Lie algebras 
$LU_\fg(\lambda)$, where $\lambda$ runs over the projective space of 
dimension equal to the rank of $\fg$, are quantizations of the Lie 
algebras of functions on the orbits of the coadjoint representation of 
$\fg$.

The Lie algebra $\fgl(\lambda)$ of matrices of complex size is the 
simplest example; it is $LU_{\fsl(2)}(\lambda)$.  The dynamical 
systems associated with it in the space of pseudodifferential 
operators in the same way as the KdV hierarchy is associated with 
$\fsl(n)$ are those studied by Gelfand--Dickey and Khesin--Malikov.  
For $\fg\neq \fsl(2)$ we get generalizations of $\fgl(\lambda)$ and 
the corresponding dynamical systems, in particular, their superized 
versions. The algebras $LU_{\fsl(2)}(\lambda)$ posess a trace and an invariant 
symmetric bilinear form, hence, with these Lie algebras associated are 
analogs of the Yang-Baxter equation, KdV, etc.

Our presentation of $LU_{\fs}(\lambda)$ for a simple $\fs$ is related to 
presentation of $\fs$ in terms of a certain pair of generators. For 
$\fs=\fsl(n)$ there are just 9 such relations.
\end{abstract}

\maketitle 

This is our paper published in: by E.~Ram\'irez de Arellano, et.  al. 
(eds.)  Proc.  Internatnl.  Symp.  Complex Analysis and related
topics, Mexico, 1996, Birkhauser Verlag, 1999, 73--105.  We just wish
to make it more accessible.  Here we made minory corrections, e.g.,
replaced $\fp\fgl(\lambda)$ with $\fsl(\lambda)$: to write
$\fp\fgl(\lambda)$ is correct, but notation $\fsl(\lambda)$ is closer
to its finite dimensional particular case.

\section*{\S 0. Introduction}

This is an expanded transcript of the talk given at the 
International Symposium on Complex Analysis and Related Topics, 
Cuernavaca, Mexica, November 18 -- 22, 1996.  We are thankful to 
A.~Turbiner and N.~Vasilevski for hospitality.

\ssec{0.0.  History} About 1966, V.~Kac and B.~Weisfeiler began the 
study of simple {\it filtered} Lie algebras of {\it polynomial 
growth}.  Kac first considered the $\Zee$-{\it graded} Lie algebras 
associated with the filtered ones and classified {\it simple graded} 
Lie algebras of {\it polynomial growth} under a technical assumption 
and conjectured the inessential nature of the assumption.  It took 
more than 20 years to get rid of the assumption: see very complicated
papers by O.~Mathieu, cf.  \cite{K} and references therein.  For a 
similar list of simple $\Zee$-graded Lie {\it super}algebras of 
polynomial growth see \cite{KS}, \cite{LSc}.

The Lie algebras Kac distinguished (or rather the algebras of 
derivations of their nontrivial central extensions, the {\it 
Kac--Moody} algebras) proved very interesting in applications.  These 
algebras aroused such interest that the study of filtered algebras was 
arrested for two decades.  Little by little, however, the simplest 
representative of the new class of simple filtered Lie superalgebras 
(of polynomial growth), namely, the Lie algebra $\fgl(\lambda)$ of 
matrices of complex size, and its projectivization, i.e., the quotient 
modulo the constants, $\fp\fgl(\lambda)$, drew its share of attention 
\cite{F}, \cite{KM}, \cite{KR}.

While we typed this paper, Shoikhet \cite{Sh} published a description 
of representations of $\fgl(\lambda)$; we are thankful to M.~Vasiliev 
who informed us of still other applications of generalizations of 
$\fgl(\lambda)$, see \cite{BWV}, \cite{KV}.

This paper begins a systematic study of a new class of Lie algebras: 
simple filtered Lie algebras of polynomial growth (SFLAPG) for which 
the graded Lie algebras associated with the filtration considered are 
not simple; $\fsl(\lambda)$ is our first example.  Actually, an
example of a Lie algebra of class SFLAPG was known even before the 
notion of Lie algebras was introduced.  Indeed, the only deformation 
(physicists call it {\it quantization}) $Q$ of the Poisson Lie algebra 
$\fpo(2n)$ sends $\fpo(2n)$ into $\fdiff(n)$, the Lie algebra of 
differential operators with polynomial coefficients; the restriction 
of $Q$ to $\fh(2n)=\fpo(2n)/center$, the Lie algebra of Hamiltonian 
vector fields, sends $\fh(2n)$ to the projectivization 
$\fp\fdiff(n)=\fdiff(n)/\Cee \cdot 1$ of $\fdiff(n)$.  The Lie algebra 
$\fp\fdiff(n)$ escaped Kac's classification, though it is the deform 
of an algebra from his list, because its intrinsically natural 
filtration given by $\deg q_i=-\deg \partial_{q_{i}}=1$ is not of 
polynomial growth while the graded Lie algebra associated with the 
filtration of polynomial growth (given by $\deg q_i=\deg 
\partial_{q_{i}}=1$) is not simple.

Observe that from the point of view of dynamical systems the Lie 
algebra $\fdiff(n)$ is not very interesting: it does not possesses a 
nondegenerate bilinear symmetric form; we will consider its 
subalgebras that do.  

\vskip 0.2 cm

In what follows we will usually denote the associative (super)algebras 
by Latin letters; the Lie (super)algebras associated with them by 
Gothic letters; e.g., $\fgl(n)=L(\Mat(n))$, $\fdiff(n)=L(\Diff(n))$, 
where the functor $L$ replaces the dot product withthe bracket.

\ssec{0.1.  The construction.  Problems related} Each of our Lie 
algebras (and Lie superalgebras) $LU_\fg(\lambda)$ is realized as a 
quotient of the Lie algebra of global sections of the sheaf of twisted 
$D$-modules on the flag variety, cf.  \cite{Ka}, \cite{Di}.  The general
construction  consists of the preparatory step 0), the main steps 1)
and 2) and two  extra steps 3) and 4).

We distinguish two cases: A) $\dim \fg<\infty$ and $\fg$ possesses a 
Cartan matrix and B) $\fg$ is a simple vectorial Lie (super)algebra.

Let $\fg=\fg_-\oplus 
\fh\oplus\fg_+$, where $\fg_+= 
\mathop{\oplus}\limits_{\alpha>0}\fg_\alpha$ and $\fg_-= 
\mathop{\oplus}\limits_{\alpha<0}\fg_\alpha$, be one of the simple 
$\Zee $-graded Lie algebras of polynomial growth, either finite 
dimensional or of vector fields, represented as the sum of its maximal 
torus (usually identical with the Cartan subalgebra) $\fh$ and the 
root subspaces $\fg_\alpha$ corresponding to an order in the set $R$ 
of roots.

Observe that each order of $R$ is in one-to-one correspondence with a 
system of simple roots.  For the finite dimensional Lie algebras $\fg$ 
all systems of simple roots are equivalent, the equivalence is 
established by the Weyl group.  For Lie superalgebras and infinite 
dimensional Lie algebras of vector fields there are inequivalent 
systems of simple roots; nevertheless, there is an analog of the Weyl 
group and the passage from system to system is described in \cite{PS}.  

For vectorial Lie algebras and Lie superalgebras, even the dimension
of the superspaces $X=(\fg_-)^*$ associated with systems of simple
roots can vary.  It is not clear if only {\it essential} (see
\cite{PS}) systems of simple roots are essential in the construction
of Verma modules (roughly speaking, each Verma module is isomorphic to
the space of functions on $X$) in which we will realize
$LU_{\fg}(\lambda)$, but hopefully not all.

\underline{Step 0): From $\fg$ to $\tilde\fg$}.  From representation
theory it is clear that there exists a realization of the elements of
$\fg$ by differential operators of degree $\leq 1$ on the space
$X=(\fg_-)^*$.  The realization has rank $\fg$ parameters (coordinates
$\lambda=(\lambda_1, \dots , \lambda_n)\in\fh^*$ of the highest weight
of the $\fg$-module $M^\lambda$).  For the algorithms of construction
and its execution in some cases see \cite{BMP}, \cite{B}, \cite{BGLS}.

Let $\tilde\fg$ be the image of $\fg$ with respect to this
realization.  Let $\tilde S^{\bcdot}(\tilde\fg)$ be the associative
subalgebra generated by $\tilde\fg$.  Clearly,
$\tilde S^{\bcdot}(\tilde\fg)\subset \fdiff(\fg_-)$.  Set
$$
U_\fg(\lambda)=\tilde 
S^{\bcdot}(\tilde\fg)/J(\lambda), \text{ where $J(\lambda)$ is the 
maxiam ideal}.
$$ 
Observe that $J(\lambda)=0$ for generic $\lambda$. 

Roughly speaking, $U_\fg(\lambda)$ is \lq\lq$\Mat"(L^\lambda)$, where
$L^\lambda$ is the quotient of $M^\lambda$ modulo the maximal
submodule $I(\lambda)$ (it can be determined and described with the
help of the Shapovalov form, see \cite{K}) and $\tilde
S^{\bcdot}(\tilde\fg)$ is the subalgebra generated by $\tilde\fg$ in
the symmetric algebra of $\tilde\fg$ modulo the relations between
differential operators.  Clearly, $\tilde S^{\bcdot}(\tilde\fg)$ is
smaller than $S^{\bcdot}(\fg)$ due to the relations between the
differential operators that span $\tilde\fg$.  

To explicitly describe the generators of $J(\lambda)$ is a main
technical problem.  We solve it in this paper for $\rk \fg=1$.  The
general case will be considered elsewhere.

\underline{Step 1) From $U_\fg(\lambda)$ to $LU_\fg(\lambda)$} Recall 
that $LU_\fg(\lambda)$ is the Lie algebra whose space is the same as 
that of $U_\fg(\lambda)$ and the bracket is the commutator.  

\underline{Step 2) Montgomery's functor} S.~Montgomery suggested
\cite{M} a construction of simple Lie superalgebras:
$$
\text{$Mo$: a central simple $\Zee$-graded algebra 
$\mapsto$ a simple Lie superalgebra.}\eqno{(Mo)}
$$
Observe that the associative algebras $U_\fg(\lambda)$ constructed
from simple Lie algebras $\fg$ are central simple.  In \cite{LM} we
intend to consider {\it Montgomery superalgebras} $Mo(U_\fg(\lambda))$
and compare them with the Lie superalgebras $LU_\fs(\lambda)$
constructed from Lie superalgebras $\fs$.  Montgomery functor often
produces new Lie superalgebras, e.g., if $\fg$ is equal to $\ff_4$ or
$\fe_i$, though not always: $Mo(U_{\fsl(2)}(\lambda))\cong
LU_{\fosp(1|2)}(\lambda)$.

\underline{Step 3) Twisted versions} An outer automorphism $a$ of
$\fG=LU_\fg(\lambda)$ or $Mo(U_\fg(\lambda))$ might single out a new
simple Lie subsuperalgebra $a_0(\fG)$, the set of fixed points of
$\fG$ under $a$.

For example, the intersection of $LU_{\fsl(2)}(\lambda)$ with the set 
of skew-adjoint differential operators is a new Lie algebra 
$\fo/\fsp(\lambda)$ while the intersection of 
$Mo(U_{\fsl(2)}(\lambda))=LU_{\fosp(1|2)}(\lambda)$ with the set of 
superskew-adjoint operators is the Lie superalgebra 
$\fo\fsp(\lambda+1|\lambda)$.  For the description of the outer 
automorphisms of $\fgl(\lambda)$ see \cite{LAS}. In general even the 
definition is unclear.

\underline{Step 4) Deformations} The deformations of Lie algebras and 
Lie superalgebras
obtained via steps 1) -- 3) may lead to new algebras of class SFLAPG,
cf.  \cite{Go}.  A.~Sergeev posed the following  interesting problšem:

{\sl what Lie algebras and Lie superalgebras can we get by applying the
above constructions 1) -- 3) to the quantum deformation $U_q(\fg)$ of
$U(\fg)$?}

\footnotesize
\begin{Remark} The above procedure can be also applied to (twisted) 
loop algebras $\fg=\fh^{(k)}$ and the stringy algebras; the result 
will be realized with differential operators of infinitely many 
indeterminates; they remind vertex operators. The algebra
$LU_{\fh^{(k)}}(\lambda)$ is a polynomial one but not of
polynomial growth. 
\end{Remark} 
\normalsize

\ssec{0.2.  Another description of $U_\fg(\lambda)$} For the finite
dimensional simple $\fg$ there is an alternative description of
$U_\fg(\lambda)$ as the quotient of $U(\fg)$ modulo the central
character, i.e., modulo the ideal $C_{(\lambda)}$ generated by rank
$\fg$ elements $C_i-k_i(\lambda)$, where the $C_i$ is the $i$-th
Casimir element and the $k_i(\lambda)$ is the (computed by
Harish-Chandra and Berezin) value of $C_i$ on $M^\lambda$.  This
description of $U_\fg(\lambda)$ goes back, perhaps, to Kostant, cf. 
\cite{Ka}.  From this description it is clear that, after the shift by
$\rho$, the half sum of positive roots, we get
$$
LU_\fg(\sigma(\lambda))\cong LU_\fg(\lambda)\quad \text{for any}\; 
\sigma\in W(\fg).
$$ 
A similar isomorphism holds for $Mo(U_\fg(\lambda))$.  In particular,
over $\Ree$, it suffices to consider the $\lambda$ that belong to one
Weyl chamber only.

For vectorial Lie algebras the description of $U_\fg(\lambda)$ as
$U(\fg)/C_{(\lambda)}$ is inapplicable.  For example, let $\fg=\fvect(n)$. 
The highest weight Verma modules are (for the standard filtration of
$\fg$) identical with Verma modules over $\fsl(n+1)$, but the center
of $U(\fvect(n))$ consists of constants only.  It is a research
problem to describe the generators of $C_{(\lambda)}$ in such cases.

Though the center of $U(\fg)$ is completely described by A.~Sergeev
for all simple finite dimensional Lie superalgebras \cite{S}, the
problem 

{\sl describe the {\it generators of the ideal} $C_{(\lambda)}$}

is open for Lie superalgebras $\fg$ even if $\fg$ is of the form
$\fg(A)$ (i.e., if $\fg$ has Cartan matrix $A$) different from
$\fosp(1|2n)$: for them the center of $U(\fg)$ is not noetherian and
it is {\it a priori} unclear if $C_{(\lambda)}$ has infinitely or
finitely many generators.  (As we will show elsewhere, $C_{(\lambda)}$
is generated for Lie superalgebras $\fg$ of the form $\fg(A)$ by the
first $\rk \fg$ Casimir operators and finitely many extra elements. 
For algebras $\fg$ of other types we do not even have a conjecture.)

\ssec{0.3.  Our result} The main result is the statement of the fact
that the above constructions 1) -- 4) yield a new class of simple Lie
(super) algebras of polynomial growth (some of which have nice
properties).

Observe that our Lie algebras $LU_{\fg}(\lambda)$ are quantizations of
the Lie algebras considered in \cite{DGS} which are also of class
SFLAPG and are contractions of our algebras.  Indeed, Donin, Gurevich
and Shnider consider the Lie algebras of functions on the orbits of
the coajoint representation of $\fg$ with respect to the Poisson
bracket.  These DGS Lie algebras are naturally realized as the
quotients of the polynomial algebra modulo an inhomogeneous ideal that
singles out the orbit; we realize the result of quantization of DGS Lie
algebras by differential operators. 

In this paper  we consider the simplest case of the superization of
this  construction: replace $\fsl(2)$ with $\fosp(1|2)$.  The cases of 
higher ranks will be considered elsewhere.  The  Khesin--Malikov
construction \cite{KM} can be applied almost literally to 
the Lie (super)algebras $LU_\fg(\lambda)$ such that $\fg$ 
admits a (super)principal embedding, see, e.g.,  \cite{GL2}.

Our main theorems: 2.6 and 4.3.  The structure of the algebras 
$LU_\fg(\lambda)$ (real forms, automorphisms, root systems) will be 
described elsewhere, see e.g., \cite{LAS}.

Observe that while the polynomial Poisson Lie algebra has only one
class of nontrivial deformations and all the deformed algebras are
isomorphic, cf.  \cite{LSc1}, the dimension of the space of parameters
of deformations of Lie algebras of Donin, Gurevich and Shnider is
equal to the rank of $\fg$ and all of the deforms are pairwise
nonisomorphic, generally.

\ssec{0.4.  The defining relations} The notion of defining relations 
is clear for a nilpotent Lie algebra.  This is one of the reasons why 
the most conventional way to present a simple Lie algebra $\fg$ is to 
split it into the direct sum of a (commutative) Cartan subalgebra and 2 
maximal nilpotent subalgebras $\fg_{\pm}$ (positive and negative).  
There are about $(2\cdot\rk\fg)^2$ relations between the $2\cdot\rk 
\fg$ generators of $\fg_{\pm}$.  The generators of $\fg_{+}$ together 
with the generators of $\fg_{-}$ generate $\fg$ as well.  In $\fg$, 
there are about $(3\cdot\rk\fg)^2$ relations between these generators; 
the relations additional to those in $\fg_{+}$ or $\fg_{-}$, i.e., 
between the positive and the negative generators, are easy to grasp.  
Though numerous, all these relations --- called {\it Serre relations} 
--- are neat and this is another reason for their popularity.  These 
relations are good to deal with not only for humans but for computers 
as well, cf.  sec.  7.3.

Nevertheless, it so happens that the Chevalley-type generators and, 
therefore, the Serre relations are not always available.  Besides, as 
we will see, there are problems in which other generators and 
relations naturally appear, cf. \cite{GL2}.

Though not so transparent as for nilpotent algebras, the notion of 
generators and relations makes sense in the general case.  For 
instance, with the principal embeddings of $\fsl (2)$ into $\fg$ one 
can associate only {\bf two} elements that generate $\fg$; we call 
them {\it Jacobson's generators}, see \cite{GL1}.  We explicitly describe 
the associated with the principal embeddings of $\fsl (2)$ 
presentations of simple Lie algebras, finite dimensional and certain 
infinite dimensional; namely, the Lie algebra \lq \lq of matrices of 
complex size" realized as a subalgebra of the Lie algebra $\fdiff(1)$ 
of differential operators in 1 indeterminate or of $\fgl_+(\infty)$, 
see \S 2.

The relations obtained are rather simple, especially for 
nonexceptional algebras.  In contradistinction with the conventional 
presentation there are just 9 relations between Jacobson's generators 
for $\fsl(\lambda)$ series (actually, 8 if 
$\lambda\in\Cee\setminus\Zee $) and not many more for the other 
algebras.

It is convenient to present $\fsl(\lambda)$ as the Lie algebra 
generated by two differential operators: $X^+=u^2\frac{d}{du}-(\lambda 
-1)u$ and $Z_{\fsl}= \frac{d^2}{du^2}$; its Lie subalgebra 
$\fo{/}\fsp(\lambda)$ of skew-adjoint operators --- a hybrid of Lie 
algebras of series $\fo$ and $\fsp$ (do not confuse with the Lie 
superalgebra of $\fosp$ type!)  --- is generated by the same $X^+$ and 
$Z_{\fo{/}\fsp}= \frac{d^3}{du^3}$; to make relations simpler, we 
always add the third generator $X^-=-\frac{d}{du}$.  For integer 
$\lambda$ each of these algebras has an ideal of finite codimension 
and the quotient modulo the ideal is the conventional $\fsl(n)$ (for 
$\lambda=n$ and $\fgl(\lambda)$) and either $\fo(2n+1)$ (for $\lambda=2n+1$) or $\fsp(2n)$ (for 
$\lambda=2n$), respectively , for $\fo/\fsp(\lambda)$.

In this paper we superize \cite{GL1}: replace $\fsl(2)$ with its 
closest relative, $\fosp(1|2)$.  We denote by 
$\fsl(\lambda|\lambda+1)$ the Lie superalgebra generated by 
$\nabla^+=x\partial_\theta +x\theta\partial _x-\lambda\theta$, 
$Z=\partial_x\partial_\theta - \theta {\partial_x}^2$ and 
$U=\partial_\theta - \theta \partial_x$, where $x$ is an even 
indeterminate and $\theta$ is an odd one.  We define 
$\fosp(\lambda+1|\lambda)$ as the Lie subsuperalgebra of 
$\fsl(\lambda|\lambda+1)$ generated by $\nabla^+$ and $Z$.  The 
presentations of $\fsl(\lambda|\lambda+1)$ and 
$\fosp(\lambda+1|\lambda)$ are associated with the {\it 
superprincipal} embeddings of $\fosp(1|2)$.  For $\lambda\in\Cee 
\setminus \Zee $ these algebras are simple.  For integer $\lambda=n$ 
each of these algebras has an ideal of finite codimension and the 
quotient modulo the ideal is the conventional $\fsl(n|n+1)$ and 
$\fosp(2n+1|2n)$, respectively.

\ssec{0.5.  Some applications} (1) Integrable systems like continuous 
Toda lattice or a generalization of the Drinfeld--Sokolov construction 
are based on the superprincipal embeddings in the same way as the 
Khesin--Malikov construction \cite{KM} is based on the principal 
embedding, cf.  \cite{GL2}.

(2) To $q$-quantize the Lie algebras of type $\fsl(\lambda)$ \`a la 
Drinfeld, using only Chevalley generators, is impossible; our 
generators indicate a way to do it.

\ssec{0.6.  Related topics} We would like to draw attention of the 
reader to several other classes of Lie algebras.  One of the reasons 
is that, though some of these classes have empty intersections with 
the class of Lie algebras we consider here, they naturally spring to 
mind and are, perhaps, deformations of our algebras in some, yet 
unknown, sense.

$\bullet$ {\it Krichever--Novikov algebras}, see \cite{SH} and refs.  
therein.  The KN-algebras are neither graded, nor filtered (at least, 
wrt the degree considered usually).  Observe that so are our algebras 
$LU_\fg(\lambda)$ with respect to the degree induced from $U(\fg)$, so 
a search for a better grading is a tempting problem.

$\bullet$ {\it Odessky} or {\it Sklyanin algebras}, see \cite{FO} and 
refs.  therein.

$\bullet$ {\it Continuum algebras}, see \cite{SV} and refs.  therein.  
In particular cases these algebras coincide with Kac--Moody or loop 
algebras, i.e., have a continuum analog of the Cartan matrix.  But to 
suspect that $\fgl(\lambda)$ has a Cartan matrix is wrong, see sec.  
2.2.  Nevertheless, in the simplest cases, if $\rk \fg=1$, the algebras 
$LU_\fg(\lambda)$ and their ``relatives'' obtained in steps 1) -- 3) 
(and, perhaps, 4)) of sec 0.1 do possess Saveliev-Vershik's nonlinear 
{\it Cartan operator} which replaces the Cartan matrix.

\section*{\S 1. Recapitulation: finite dimensional simple Lie algebras}

This section is a continuation of \cite{LP}, where the case of the 
simplest base (system of simple roots) is considered and where 
non-Serre relations for simple Lie algebras first appear, though in a 
different setting.  This paper is also the direct superization of 
\cite{GL1}; we recall its results.  For presentations of Lie superalgebras 
with Cartan matrix via Chevalley generators, see \cite{LS}, \cite{GL3}.
 
What are ``natural'' generators and relations for a {\it simple finite
dimensional Lie algebra}?  The answer is important in questions when
it is needed to identify an algebra $\fg$ given its generators and
relations.  (Examples of such problems are connected with
Estabrook--Wahlquist prolongations, Drinfeld's quantum algebras,
symmetries of differential equations, integrable systems, etc.).
  
\ssec{1.0.  Defining relations} If $\fg$ is nilpotent, the problem of 
its presentation has a natural and unambiguous solution: 
representatives of the homology $H_1(\fg)\cong \fg/ [\fg, \fg]$ are 
the generators of $\fg$ and the elements from $H_2(\fg)$ correspond to 
relations.

On the other hand, if $\fg$ is simple, then $\fg =[\fg, \fg]$ and 
there is no ``most natural'' way to select generators of $\fg$.  The 
choice of generators is not unique.

Still, among algebras with the property $\fg =[\fg, \fg]$ the simple 
ones are distinguished by the fact that their structure is very well 
known.  By trial and error people discovered that for finite 
dimensional simple Lie algebras, there are certain ``first among 
equal'' sets of generators:

1) {\it Chevalley generators} corresponding to positive and negative 
simple roots;

2) a pair of generators that generate any finite dimensional simple 
Lie algebra associated with the {\it principal $\fsl(2)$-subalgebra} 
(considered below).

The relations associated with Chevalley generators are well-known, see 
e.g., \cite{OV}, \cite{K}.  These relations are called {\it Serre relations}.

The possibility to generate any simple finite dimensional Lie algebra 
by two elements was first claimed by N.~Jacobson; for the first (as 
far as we know) proof see \cite{BO}.  We do not know what generators 
Jacobson had in mind; \cite{BO} take for them linear combinations of 
positive and negative root vectors with generic coefficients; nothing 
like a \lq\lq natural" choice that we suggest to refer to as {\it 
Jacobson's generators} was ever proposed.

To generate a simple algebra with only two elements is tempting but 
nobody yet had explicitly described relations between such generators, 
perhaps, because to check whether the relations between these elements 
are nice-looking is impossible without a modern computer (cf.  an 
implicit description in \cite{F}).  As far as we could test, the relations 
for any other pair of generators chosen in a way distinct from ours 
are too complicated.  There seem to be, however, one exception cf. 
\cite{GL2}.

\ssec{1.1.  The principal embeddings} There exists only one (up to 
equivalence) embedding $r: \fsl(2)\tto \fg$ such that $\fg$, 
considered as $\fsl(2)$-module, splits into $\rk \fg$ irreducible 
modules, cf.  \cite{D} or \cite{OV}.  This embedding is called {\it 
principal} and, sometimes, {\it minimal} because for the other 
embeddings (there are plenty of them) the number of irreducible 
$\fsl(2)$-modules is $>\rk \fg$.  Example: for $\fg=\fsl(n)$, 
$\fsp(2n)$ or $\fo(2n+1)$ the principal embedding is the one 
corresponding to the irreducible representation of $\fsl(2)$ of 
dimension $n$, $2n$, $2n+1$, respectively.

For completeness, let us recall how the irreducible $\fsl(2)$-modules
with highest weight look like.  (They are all of the form $L^{\mu}$,
where $L^{\mu}=M^{\mu}$ if $\mu\not\in \Zee _+$, and
$L^n=M^{n}/M^{-n-2}$ if $n\in \Zee _+$, and where $M^{\mu}$ is
described below.)  Select the following basis in $\fsl(2)$:
$$ 
X^-=\pmatrix
0&0\\ -1&0\cr\endpmatrix, \quad H=\pmatrix 1&0\\  0&-1\\ \endpmatrix, \quad X^+= \pmatrix
0&1\\  0&0\\ 
\endpmatrix.
$$
The $\fsl(2)$-module $M^{\mu}$ is illustrated with a graph whose nodes 
correspond to the eigenvectors $l_{\mu-2i}$ of $H$ with the weight 
indicated;
$$ 
\dots\overset{\mu-2i-2}{\circ} -\overset{\mu-2i}{\circ} -\dots 
-\overset{\mu-2}{\circ} -\overset{\mu}{ \circ}
$$
the edges depict the action of $X^{\pm}$ (the action of
$X^+$ is directed to the right, that of $X^-$ to the left: $X^-l_{\mu-2i}=l_{\mu-2i-2}$ and
$$
X^+l_{\mu-2i}=X^+((X^-)^il_{\mu})=i(\mu-i+1)l_{\mu-2i+2};\quad
X^+(l_{\mu})=0.\eqno{(1.1)}
$$   
 
As follows from (1.1), the module $M^n$ for $n\in \Zee _+$ has an
irreducible submodule isomorphic to $M^{-n-2}$; the quotient,
obviously irreducible, as follows from the same (1.1), will be denoted
by $L^n$.

There are principal $\fsl(2)$-subalgebras in every finite dimensional
simple Lie algebra, though, generally, not in infinite dimensional
ones, e.g., not in affine Kac-Moody algebras.  The construction is as
follows.  Let $X^{\pm}_1, \dots , X^{\pm}_{\rk \fg}$ be Chevalley
generators of $\fg$, i.e., the generators corresponding to simple
roots.  Let the images of $X^{\pm}\in\fsl(2)$ in $\fg$ be
$$
X^-\mapsto \sum X^{-}_i;\quad X^+\mapsto \sum a_iX^{+}_i
$$
and select the $a_i$ from the relations $[[X^+, X^-], X^{\pm}]=\pm
2X^{\pm}$ true in $\fsl(2)$.  For $\fg$ constructed from a Cartan
matrix $A$, there is a solution for $a_i$ if and only if $A$ is invertible.

In Table 1.1 a simple finite dimensional Lie algebra $\fg$ is 
described as the $\fsl(2)$-module corresponding to the principal 
embedding (cf.  \cite{OV}, Table 4).  The table introduces the number 
$2k_2$ used in relations.  We set $k_{1}=1$.

\ssec{Table 1.1. $\fg$ as the $\fsl(2)$-module}
$$
\renewcommand{\arraystretch}{1.4}
\begin{tabular}{|l|l|r|}
\hline
$\fg$&the $\fsl(2)$-spectrum of $\fg=L^2\oplus L^{2k_2} 
\oplus L^{2k_3}
\dots$ &$2k_2$\cr
\hline
$\fsl(n)$&$L^2\oplus L^4 \oplus L^6 \dots \oplus L^{2n-2}$&4\cr
$\fo (2n+1)$, $\fsp(2n)$\; &$L^2\oplus L^6 \oplus L^{10} \dots \oplus L^{4n-2}$&6\cr
$\fo (2n)$&$L^2\oplus L^6 \oplus L^{10} \dots \oplus L^{4n-6}\quad\oplus L^{2n-2}$&6\cr
$\fg_2$&$L^2\oplus L^{10}$&10\cr
$\ff_4$&$L^2\oplus L^{10}\oplus L^{14}\oplus L^{22}$&10\cr
$\fe_6$&$L^2\oplus L^{8}\oplus L^{10}\oplus L^{14}\oplus L^{16}\oplus L^{22}$&8\cr 
$\fe_7$&$L^2\oplus L^{10}\oplus L^{14}\oplus L^{18}\oplus L^{22}\oplus L^{26}\oplus
L^{34}$&10\cr 
$\fe_8$&$L^2\oplus L^{14}\oplus L^{22}\oplus L^{26}\oplus L^{34}\oplus L^{38}
\oplus L^{46}\oplus L^{58}$&14\cr
\hline
\end{tabular}
$$

One can show that $\fg$ can be generated by two elements: $x:=X^+\in
L^2=\fsl (2)$ and a lowest weight vector $z:=l_{-r}$ from an
appropriate module $L^r$ other than $L^2$ from Table 1.1.  For the
role of this $L^r$ we take either $L^{2k_{2}}$ if $\fg\not=\fo(2n)$ or
the last module $L^{2n-2}$ in the above table if $\fg =\fo(2n)$. 
(Clearly, $z$ is defined up to proportionality; we will assume that a
basis of $L^r$ is fixed and denote $z=t\cdot l_{-r}$ for some $t\in
\Cee $ that can be fixed at will, cf.  \S 3.)

The exceptional choice for $\fo(2n)$ is occasioned by the fact that by
choosing $z\in L^{r}$ for $r\neq 2n-2$ instead, we generate
$\fo(2n-1)$.

We call the above $x$ and $z$, together with $y:=X^-\in L^2$ taken for
good measure, {\it Jacobson's generators}.  The presence of $y$
considerably simplifies the form of the relations, though slightly
increases their number.  (One might think that taking the symmetric to
$z$ element $l_r$ will improve the relations even more but in reality
just the opposite happens.)

Concerning $\fg =\fo(2n)$ see sec. 7.2.

\ssec{1.2.  Relations between Jacobson's generators} First, observe
that if an ideal of a free Lie algebra is homogeneous (with respect to
the degrees of the generators of the algebra), then the number and the
degrees of the defining relations (i.e., the generators of the ideal)
is uniquely defined provided the relations are homogeneous.  This is
obvious.

A simple Lie algebra $\fg$, however, is the quotient of a free Lie
algebra $\fF$ modulo a inhomogeneous ideal, $\fI$, the ideal without
homogeneous generators.  Therefore, we can speak about the number and
the degrees of relations only conditionally.  Our condition is the
possibility to express any element $x\in\fI$ via the generators $g_1,
...$ of $\fI$ by a formula of the form
$$
x=\sum [c_i,g_i],\text{ where $c_i\in \fF$ and $\deg c_i+\deg g_i\leq 
\deg x$ for all $i$.} \eqno{(*)}
$$

Under condition $(*)$ the number of relations and their degrees are 
uniquely determined.  Now we can explain why do we need an extra 
generator $y$: without $y$ the weight relations would have been of 
very high degree.

We divide the relations between the Jacobson generators into the types
corresponding to the number of occurrences of $z$ in them: {\bf 0}. 
Relations in $L^2 = \fsl(2)$; {\bf 1}.  Relations coming from the
$\fsl(2)$-action on $L^{2k_2}$; {\bf 2}.  Relations coming from
$L^{2k_{1}}\wedge L^{2k_2}$; $\pmb{\geq 3}$.  Relations coming from
$L^{2k_2}\wedge L^{2k_2}\wedge L^{2k_2}\wedge \dots$ with $\geq 3$
factors; among the latter relations we distinguish one --- of type
``$\pmb\infty$" --- the relation that shears the dimension.  (For small
$\text{rank}\, \fg$ the relation of type $\infty$ can be of the above
types.)

Observe that, apart form relations of type $\infty$, the relations of
type $\geq 3$ are those of type 3 except for $\fe_7$ which satisfies
stray relations of types 4 and 5, cf.  \cite{GL1}.

The relations of type 0 are the well-known relations in $\fsl(2)$
$$
\pmb{0.1}. \; \; [[x, y], \, x]=2x,\quad \quad  \quad   
\pmb{0.2}\; \;  [[x, y], \, y]=-2y.\eqno{(Rel 0)}
$$
The relations of type 1 mirror the fact that the space $L^{2k_2}$ is 
the $(2k_2+1)$-dimensional $\fsl(2)$-module.  To simplify notations we 
denote: $z_i=(\ad x)^iz$.  Then the type 1 relations are:

$$
\pmb{1.1}.\; \;  [y, \, z] = 0,\; \pmb{1.2}. \; \;  [[x, y], \, z] = 
-2k_2z,\; 
\pmb{1.3}.\; \;  z_{2k_{1}+1} = 0\; \text{ with $2k_2$ from Table
1.1.}.\eqno{(Rel 1)} 
$$

\ssbegin{1.3}{Theorem} For the simple finite dimensional Lie algebras
all the relations between the Jacobson generators are the above
relations {\em (Rel0), (Rel1)} and the relations from {\em
\cite{GL1}}.  \end{Theorem} 

In \S 3 these relations from \cite{GL1} are reproduced for
the classical Lie algebras.

\section*{\S 2.  The Lie algebra $\fsl(\lambda)$ as a quotient 
algebra of $\fdiff (1)$ and a subalgebra of $\fsl_+ (\infty)$}

\ssec{2.1.  $\fgl(\lambda)$ is endowed with a trace} The Poincar\' 
e-Birkhoff-Witt theorem states that, as spaces, $U(\fsl(2))\cong \Cee 
[X^-, H, X^+]$.  We also know that to study representations of $\fg$ 
is the same as to study representations of $U(\fg)$.  Still, if we are 
interested in irreducible representations, we do not need the whole of 
$U(\fg)$ and can do with a smaller algebra, easier to study.

This observation is used now and again; Feigin applied it in \cite{F} 
writing, actually, (as deciphered in \cite{PH}, \cite{GL1}, \cite{Sh}) 
that setting
$$ 
X^-=-\frac{d}{d u}, \quad
H=2u\frac{d}{d u}-(\lambda-1),  \quad X^+=u^2\frac{d}{d
u}-(\lambda-1) u  \eqno{(2.1)}
$$  
we obtain a morphism of $\fsl(2)$-modules and, moreover, of
associative algebras: $U(\fsl(2))\longrightarrow \Cee [u, \frac{d}{d
u}]$.  The kernel of this morphism is the ideal generated by
$\Delta-\lambda ^2+1$, where $\Delta=2(X^+X^-+X^-X^+)+H^2$.  Observe,
that this morphism is not an epimorphism, either.  The image of this
morphism is our Lie algebra of matrices of \lq\lq complex size".

\begin{rem*}{Remark}
In their proof of certain statements from \cite{F} that we will 
recall, \cite{PH} made use of the well-known fact that the Casimir 
operator $\Delta$ acts on the irreducible $\fsl(2)$-module $L^{\mu}$ 
(see sec 1.1) as the scalar operator of multiplication by 
$\mu^2+2\mu$.  The passage from \cite{PH}'s $\lambda$ to \cite{F}'s 
$\mu$ is done with the help of a shift by the weight $\rho$, a half 
sum of positive roots, which for $\fsl(2)$ can be identified with 1, 
i.e., $(\lambda-1)^2+2(\lambda-1)=\lambda^2-1$ for $\lambda=\mu+1$.
\end{rem*}

Consider the Lie algebra $LU(\fsl(2))$ associated with the associative
algebra $U(\fsl(2))$.  Set
$$
U_\lambda=U(\fsl(2))/(\Delta-\lambda ^2+1). \eqno{(2.2)}
$$
The definition directly implies that
$\fgl(-\lambda)\cong\fgl(\lambda)$, so speaking about real values of
$\lambda$ we can confine ourselves to the nonnegative values, cf. 
sec.  0.2.  It is easy to see that, as $\fsl(2)$-module,
$$ 
LU_\lambda=L^0\oplus L^2\oplus L^4\oplus\dots\oplus
L^{2n}\oplus \dots \eqno{(2.3)}
$$
It is not difficult to show (see \cite{PH} for details) that the Lie
algebra $LU_n$ for $n\in \Zee \setminus\{0\}$ contains an ideal
$J_{n}$ and the quotient $LU_n/J_{n}$ is the conventional $\fgl(n)$. 
In \cite{PH} it is proved that for $\lambda \neq \Zee\setminus\{0\}$
the Lie algebra $LU_\lambda$ has only two ideals --- the space $L^0$
of constants and its complement.  Set
$$
\fp\fgl(\lambda)=\fgl(\lambda)/L^0, \text{ where }\fgl(\lambda)=
\cases LU_\lambda&\text{for $\lambda \not\in \Zee\setminus\{0\}$}\cr
LU_n/J_{n}&\text{for $n\in\Zee\setminus\{0\}$.}\endcases
\eqno{(2.4)}
$$

Observe, that $\fgl(\lambda)$ is endowed with a trace. This 
follows directly from (2.3) and the fact that 
$$
\fgl(\lambda)\cong L^0\oplus [\fgl(\lambda), \fgl(\lambda)].
$$
Therefore, $\fp\fgl(\lambda)$ can be identified with $\fsl(\lambda)$,
the subalgebra of the traceless matrices in $\fgl(\lambda)$.  We can
normalize the trace at will, for example, if we set $\tr
(id)=\lambda$, then the trace that our trace induces on the quotient
of $LU_{\fsl(2)}(n)$ modulo $J(n)$ coincides with the usual trace on
$\fgl(n)$ for $n\in\Nee$.

Another way to introduce the trace was suggested by J. Bernstein.  We
decipher its description in \cite{KM} as follows.  Look at the image
of $H\in\fsl(2)$ in $\fgl(M^\lambda)$.  Bernstein observed that though
the trace of the image is an infinite sum, the sum of the first $D+1$
summands is a polynomial in $D$, call it $\tr (H)$.  It is easy to see
that $\tr (H)$ vanishes if $D=\lambda$.

Similarly, for {\it any} $x\in LU_{\fg}(\lambda)$ considered as an 
element of $\fgl(M^\lambda)$ set 
$$
\tr(x; D)=\mathop{\sum}\limits_{i=1}^{D}x_{ii}.
$$
Let $D(\lambda)$ be the value of the dimension of the irreducible 
finite dimensional $\fg$-module with highest weight $\lambda$, for an 
exact formula see \cite{D}, \cite{OV}.  Set $\tr(x)=\tr(x; 
D(\lambda))$; as is easy to see, this formula determines the trace on 
$LU_{\fg}(\lambda)$ for arbitrary values of $\lambda$.

Observe that whereas for any irreducible finite dimensional module 
over the simple Lie algebra $\fg$ there is just one  formula for 
$D(\lambda)$ (H.~Weyl dimension formulas) for Lie superalgebra there 
are {\it several} distinct formulas depending on how ``typical'' 
$\lambda$ is.

\ssec{2.2.  There is no Cartan matrix for $\fsl(\lambda)$.  What
replaces it?} Are there {\it Chevalley generators} in $\fsl(\lambda)$? 
In other words are there elements $X^{\pm}_i$ of degree $\pm 2$ and
$H_i$ of degree 0 (the {\it degree} is the weight with respect to the
$\fsl(2)=L^2\subset \fsl(\lambda)$) such that
$$
[X^+_i, X^-_j]= \delta_{ij}H_i, \quad [H_i, H_j]=0\text{ and }[H_i,
X^{\pm}_j]=\pm A_{ij}X^{\pm}_j?  \eqno{(2.5)}
$$
The answer is {\bf NO}: $\fsl(\lambda)$ is too small.  To see what is
the problem, consider the following elements of degree $\pm 2$ from
$L^4$ and $L^6$ of $\fgl(\lambda)$:
$$
\renewcommand{\arraystretch}{1.4}
\begin{array}{ll}
\deg= -2: &-4uD^2-2(\lambda-2)D\\
\deg=2: &-4u^3D^2+6(\lambda-2)u^2D-2(\lambda-1)(\lambda-2)u 
\end{array}
$$$$
\renewcommand{\arraystretch}{1.4}
\begin{array}{ll}\deg=-2: &15u^2D^3-15(\lambda-3)uD^2+3(\lambda-2)(\lambda-3)D\\
\deg=2: &15u^4D^3-30(\lambda-3)u^3D^2+ 
18(\lambda-2)(\lambda-3)u^2D-3(\lambda-1)(\lambda-2)(\lambda-3)u
\end{array}
$$
To satisfy $(2.5)$, we can complete $\fgl(\lambda)$ by considering 
infinite sums of its elements, but the completion erases the 
difference between different $\lambda$'s:

\begin{Proposition} For $\lambda\neq \rho$ the completion of 
$\fsl(\lambda)$ generated by Jacobson's generators {\em (see Tables)} 
is isomorphic to $\overline{\fp\fdiff (1)}$, the quotient of the Lie 
algebra of differential operators with formal coefficients modulo 
constants.
\end{Proposition}

Though there is no Cartan matrix, Saveliev and Vershik \cite{SV} 
suggested an operator $K$ which replaces Cartan matrix. For furtehr 
details see paper by Shoihet and Vershik \cite{ShV}.

\ssec{2.3. The outer automorphism of $LU_{\fg}(\lambda)$} 
The invariants of the mapping 
$$
X\mapsto -SX^{t}S\;  \text{for}\;  X\in\fgl(n), \; \text{where}\;
S=\antidiag(1, -1, 1, -1\dots )\eqno{(2.6)}
$$
constitute $\fo(n)$ if $n\in 2\Nee +1$ and $\fsp(n)$ if $n\in 2\Nee $.
By analogy, Feigin defined $\fo(\lambda)$ and $\fsp(\lambda)$ as
subalgebras of $\fgl(\lambda)=\mathop{\oplus}\limits_{k\geq 0} L^{2k}$ invariant with respect to the
involution analogous to (2.6):
$$
X\mapsto \cases -X&\text{if $X\in L^{4k}$}\cr
X&\text{if $X\in L^{4k+2}$.}\endcases\eqno{(2.7)}
$$

Since $\fo(\lambda)$ and $\fsp(\lambda)$ --- the subalgebras of
$\fgl(\lambda)$ singled out by the involution (2.7) --- differ by a
shift of the parameter $\lambda$, it is natural to denote them
uniformly (but so as not to confuse with the Lie superalgebras of
series $\fosp$), namely, by $\fo{/}\fsp(\lambda)$.  For integer values
of the parameter it is clear that
$$
\fo{/}\fsp(\lambda) = \left\{\begin{matrix} \fo(\lambda)\supplus
I_\lambda& \text{if $\lambda \in 2\Nee+1$}, \\
\fsp(\lambda)\supplus I_\lambda & \text{if $\lambda \in 2\Nee
$},\end{matrix}\right.\; \text{where}\; I_\lambda \; \text{is an
ideal}.
$$
In the realization of $\fsl(\lambda)$ by differential operators the
transposition is the passage to the adjoint operator; hence,
$\fo{/}\fsp(\lambda)$ is a subalgebra of $\fsl(\lambda)$ consisting of
self-skew-adjoint operators with respect to the involution
$$
a(u)\frac{d^k}{du^k}\mapsto (-1)^k\frac{d^k}{du^k}a(u)^*.\eqno{(2.8)}
$$
The superization of this formula is straightforward: via Sign Rule.

\ssec{2.4.  The Lie algebra $\fgl(\lambda)$ as a subalgebra of
$\fgl_+(\infty)$} Recall that $\fgl_+(\infty)$ often denotes the Lie
algebra of infinite (in one direction; index $+$ indicates that)
matrices with nonzero elements inside a (depending on the matrix)
strip along the main diagonal and containing it.  The subalgebras
$\fo(\infty)$ and $\fsp(\infty)$ are naturally defined, while
$\fsl(\infty)$ is, by abuse of language, sometimes used to denote
$\fp\fgl(\infty)$. 

When it comes to superization, one shall be very careful selecting an
appropriate candidate for $\fsl(\infty|\infty)$ and its subalgebra,
cf.  \cite{E}.

The realization (2.1) provides with an embedding
$\fsl(\lambda)\subset\fsl_+(\infty)=\lq\lq\fsl(M^{\lambda})"$, so for
$\lambda\neq \Nee$ the Verma module $M^{\lambda}$ with highest weight
$\mu$ is an irreducible $\fsl(\lambda)$-module.

\begin{Proposition} The completion of $\fgl(\lambda)$ (generated by
the elements of degree $\pm 2$ with respect to $H\in
\fsl(2)\subset\fgl(\lambda)$) is isomorphic for any noninteger
$\lambda$ to $\fgl_+(\infty)=\lq\lq\fgl(M^{\lambda})"$.
\end{Proposition}

\ssec{2.5. The Lie algebras $\fsl(*)$ and $\fo{/}\fsp(*)$ for
$*\in\Cee P^1=\Cee \cup\{*\}$} 
The \lq\lq dequantization" of the relations for $\fsl(\lambda)$ and
$\fo{/}\fsp(\lambda)$ (see \S 3) is performed by passage to the limit as
$\lambda\longrightarrow\infty$ under the change:  
$$ 
t\mapsto
\left\{\renewcommand{\arraystretch}{1.4}
\begin{array}{ll}
\frac{t}{\lambda}&\text{for $\fsl(\lambda)$}\\
\frac{t}{\lambda^2}&\text{for $\fo{/}\fsp(\lambda)$}.\end{array}\right.
$$ 
So the parameter $\lambda$ above can actually run over $\Cee 
P^1=\Cee\cup\{*\}$, not just $\Cee$.  In the realization with the help 
of deformation, cf.  2.7 below, this is obvious.  Denote the limit 
algebras by $\fsl(*)$ and $\fo{/}\fsp(*)$ in order to distinguish 
them from $\fsl(\infty)$ and $\fo(\infty)$ or $\fsp(\infty)$ from sec.  
2.4.

It is clear that it is impossible to embed  $\fsl(*)$ and $\fo{/}\fsp(*)$ into the \lq\lq
quadrant" algebra $\fsl_+(\infty)$: indeed, $\fsl(*)$ and
$\fo{/}\fsp(*)$ are subalgebras of the whole
\lq\lq plane" algebras $\fsl(\infty)$ and
$\fo(\infty)$ or $\fsp(\infty)$.

\ssbegin{2.6}{Theorem} For Lie algebras $\fsl(\lambda)$ and
$\fo{/}\fsp(\lambda)$, $\lambda\in \Cee P^1$, all the relations
between the Jacobson generators are the relations of types $0, 1$ with
$2k_2$ found from Table $1.1$ and the borrowed from {\em [GL1]}
relations from \S $3$.
\end{Theorem}

\section*{\S 3. Jacobson's generators and relations between them}

In what follows the $E_{ij}$ are the matrix units; $X^{\pm}_i$ stand
for the conventional Chavalley generators of $\fg$.  For
$\fsl(\lambda)$ and $\fo{/}\fsp(\lambda)$ the generators $x =
u^2\frac{d}{du} - (\lambda-1) u$ and $y = - \frac{d}{d u}$ are the
same; $z_{\fsl} = t\frac{d^2}{d u^2}$ while $z_{\fo{/}\fsp} =
t\frac{d^3}{d u^3}$.  For $n\in\Cee \setminus \Zee $ there is no
shearing relation of type $\infty$; for $n=*\in\Cee P^1$ the relations
are obtained with the substitution 2.5.  The parameter $t$ can be
taken equal to 1; we kept it explicit to clarify how to ``dequantize"
the relations as $\lambda\tto\infty$.

\noindent $\underline{\fsl(*)}$.

$$
\renewcommand{\arraystretch}{1.3}
\begin{tabular}{rl}
{\pmb 2.1}. & $3[z_1, z_2] - 2[z, z_3] = 24 y$,\cr
{\pmb 3.1}.& $[z, \, [z, z_1] ]= 0$,\cr
{\pmb 3.2}.& $4 [[z, z_1], z_3]]] + 3 [z_2, [z, z_2]]
	= -576 z$.\cr
\end{tabular}
$$

\noindent $\underline{\fo{/}\fsp(*)}$.

$$
\renewcommand{\arraystretch}{1.3}
\begin{tabular}{rl}
{\pmb 2.1}.& $2 [z_1, z_2] - [z, z_3] = 72 z$,\cr
{\pmb 2.2}.& $9 [z_2, z_3] - 5 [z_1, z_4] =
	216 z_2 - 432 y$,\cr
{\pmb 3.1}.& $[z, \, [z, z_1] ]= 0$,\cr
{\pmb 3.2}.& $7 [[z, z_1], z_3] + 6 [z_2, [z, z_2]] =
	- 720 [z, z_1] $.\cr
\end{tabular}
$$

\vskip 1mm
\noindent $\underline{\fsl(n)\text{ for }n\geq 3}$. Generators:
$$ 
x = \sum\limits_{1\leq i\leq n-1}i(n-i)E_{i, i+1}, \qquad
y = \sum\limits_{1\leq i\leq n-1}E_{i+1, i}, \qquad
z = t\sum\limits_{1\leq i\leq n-2} E_{i+2, i}.
$$

Relations:
$$
\renewcommand{\arraystretch}{1.3}
\begin{tabular}{rl}
{\pmb 2.1}.& $3[z_1, \, z_2] - 2[z, \, z_3] = 24t^2(n^2-4) y$,\cr
{\pmb 3.1}.& $[z, \, [z, z_1]] = 0$,\cr
{\pmb 3.2}.& $4 [z_3, \, [z, z_1] ] - 3 [z_2, \, [z, z_2]]
	 = 576 t^2(n^2-9)z$.\cr
${\pmb \infty=n-1}$.& $(\ad z_1)^{n-2}z = 0$.\cr
\end{tabular}
$$

\noindent For $n=3, 4$ the degree of the last relation is lower than the
degree of some other relations, this yields simplifications. 

\noindent $\underline{\fo(2n+1)\text{ for }n\geq 3}$. Generators: 
$$
x = n(n+1)(E_{n+1, 2n+1}-E_{n, n+1}) +
  \sum\limits_{1\leq i\leq n-1} i(2n+1-i)(E_{i, i+1}-E_{n+i+2, n+i+1}),  
$$
$$
y = (E_{2n+1, n+1}-E_{n+1, n}) +
  \sum\limits_{1\leq i\leq n-1}(E_{i+1, i}-E_{n+i+1, n+i+2}),  
$$
$$
z=t\bigr( (E_{2n-1, n+1}-E_{n+1, n-2})-(E_{2n+1, n-1}-E_{2n, n}) +
  \sum\limits_{1\leq i\leq n-3} (E_{i+3, i}-E_{n+i+1, n+i+4}) \bigl).
$$

Relations:
$$
\begin{tabular}{rl}
{\pmb 2.1}.& $2 [z_1, \, z_2] - [z, \, z_3] = 144 t (2n^2+2n-9) z$,\cr
{\pmb 2.2}.& $9 [z_2, \, z_3] - 5 [z_1, \, z_4] =
 	432 t (2n^2+2n-9)z_2 + 1728 t^2(n-1)(n+2)(2n-1)(2n+3) y$,\cr
{\pmb 3.1}.& $[z, \, [z, z_1]] = 0$,\cr
{\pmb 3.2}.& $7 [z_3, \, [z, z_1]] - 6 [z_2, \, [z, z_2]] =
  2880 t (n-3)(n+4)[z, z_1] $,\cr
$\pmb{\infty=n}$.& $(\ad z_1)^{n-1}z = 0$.\cr
\end{tabular}
$$

\noindent $\underline{\fsp(2n)\text{ for }n\geq 3}$. Generators: 
$$
x=n^2E_{n, 2n} +
 \sum\limits_{1\leq i\leq n-1} i(2n-i)(E_{i, i+1}-E_{n+i+1, n+i}),  
$$
$$
y=E_{2n, n} + \sum\limits_{1\leq i\leq n-1}(E_{i+1, i}-E_{n+i, n+i+1}),  
$$
$$
z=t\biggl((E_{2n, n-2}+E_{2n-2, n}) - E_{2n-1, n-1} +
 \sum\limits_{1\leq i\leq n-3} (E_{i+3, i}-E_{n+i, n+i+3}) \biggr).
$$

Relations:
$$
\begin{tabular}{rl}
{\pmb 2.1}.& $2 [z_1, \, z_2] - [z, \, z_3] = 72 t (4n^2-19) z$,\cr
{\pmb 2.2}.& $9 [z_2, \, z_3] - 5 [z_1, \, z_4] =
	216 t (4n^2-19)z_2 + 1728 t^2(n^2-1)(4n^2-9) y$,\cr
{\pmb 3.1}.& $[z, \, [z, z_1]] = 0$,\cr
{\pmb 3.2}.& $7 [z_3, \, [z, z_1] ] - 6 [z_2, \, [z, z_2]] =
	 720 t (4n^2-49)[z, z_1] $, \cr
$\pmb{\infty=n}$.& $(\ad z_1)^{n-1}z = 0$.\cr
\end{tabular}
$$

For Jacobson generators and corresponding defining relations for the 
exceptional Lie algebras see \cite{GL1}.

\section*{\S 4. Lie superalgebras}
\ssec{4.0. Linear algebra in superspaces}
Superization has certain subtleties, often disregarded or expressed
too briefly. We will dwell on them a bit, see \cite{L2}.

A {\it superspace} is a $\Zee /2$-graded space; for a superspace 
$V=V_{\bar 0}\oplus V_{\bar 1}$ denote by $\Pi (V)$ another copy of 
the same superspace: with the shifted parity, i.e., $(\Pi(V))_{\bar 
i}= V_{\bar i+\bar 1}$.

A superspace structure in $V$ induces that in the space $\End (V)$.  A 
{\it basis of a superspace} is always a basis consisting of {\it 
homogeneous} vectors; let $\Par=(p_1, \dots, p_{\dim V})$ be an 
ordered collection of their parities, called the {\it format} of (the 
basis of) $V$.  A square {\it supermatrix} of format (size) $\Par$ is 
a $\dim V\times \dim V$ matrix whose $i$th row and $i$th column are 
said to be of parity $p_i$.  The matrix unit $E_{ij}$ is supposed to 
be of parity $p_i+p_j$ and the bracket of supermatrices (of the same 
format) is defined via Sign Rule: {\it if something of parity $p$ 
moves past something of parity $q$ the sign $(-1)^{pq}$ accrues; the 
formulas defined on homogeneous elements are extended to arbitrary 
ones via linearity}.  For example: $[X, Y]=XY-(-1)^{p(X)p(Y)}YX$; the 
sign $\wedge$ in what follows is also understood in supersence, etc.

Usually, $\Par$ is considered to be of the form $(\bar 0, \dots, \bar 
0, \bar 1, \dots, \bar 1)$.  Such a format is called {\it standard}.  
The Lie superalgebra of supermatrices of size $\Par$ is denoted by 
$\fgl(\Par)$, usually $\fgl(\bar 0, \dots, \bar 0, \bar 1, \dots, \bar 
1)$ is abbreviated to $\fgl(\dim V_{\bar 0}|\dim V_{\bar 1})$.

For $\dim V_{\bar 0}=\dim V_{\bar 1}\pm 1$ we will often use another 
format, the {\it alternating} one, $\Par_{alt}=(\bar 0, \bar 1, \bar 
0, \bar 1, \dots )$.

The {\it supertrace} is the map $\fgl (\Par)\longrightarrow \Cee $, 
$(A_{ij})\mapsto \sum (-1)^{p_{i}}A_{ii}$.  The supertraceless 
matrices constitute a Lie subsuperalgebra, $\fsl(\Par)$.

To the linear map $F$ of superspaces there corresponds the dual map 
$F^*$ between the dual superspaces; if $A$ is the supermatrix 
corresponding to $F$ in a format $\Par$, then to $F^*$ the {\it 
supertransposed} matrix $A^{st}$ corresponds:
$$
(A^{st})_{ij}=(-1)^{(p_{i}+p_{j})(p_{i}+p(A))}A_{ji}.
$$

The supermatrices $X\in\fgl(\Par)$ such that 
$$
X^{st}B+(-1)^{p(X)p(B)}BX=0\quad \text{for a homogeneous matrix 
$B\in\fgl(\Par)$}
$$
constitute the Lie superalgebra $\faut (B)$ that preserves the 
bilinear form on $V$ with matrix $B$.

The superspace of bilinear forms is denoted by $Bil_C(M, N)$ or 
$Bil_C(M)$ if $M$=$N$.  The {\it upsetting of forms} $uf:\ Bil_C(M, 
N)\rightarrow Bil_C(N, M)$, is defined by the formula
$$
B^{uf}(n, m)=(-1)^{p(n) p(m) }B(m, n).
$$ 
A form $B\in Bil_C(M)$ is called {\it supersymmetric} if $B^{uf}=B$ 
and {\it superskew-symmetric} if $B^{uf}=-B$.

Given bases $\{m_i\}$ and $\{n_j\}$ of $C$-modules $M$ and $N$ and a 
bilinear form $B: M\otimes N\rightarrow C$, we assign to $B$ the 
matrix
$$
({}^{mf}\! B)_{ij}=(-1)^{p(m_i)p(B)}B(m_i, n_j).
$$

For a nondegenerate supersymmetric form whose matrix in the standard 
format is
$$
B_{m, 2n}= \begin{pmatrix}
1_m&0\\
0&J_{2n}
\end{pmatrix},\quad \text{where $J_{2n}=\begin{pmatrix}
 0&1_n\\-1_n&0\end{pmatrix}$}.
$$
The usual notation for $\faut (B_{m|2n})$ is $\fosp^{sy}(m|2n)$ or 
just $\fosp(m|2n)$.  (Observe that the passage from $V$ to $\Pi (V)$ 
sends the supersymmetric forms to superskew-symmetric ones, preserved 
by $\fosp^{sk}(m|2n)$ which is isomorphic to $\fosp(m|2n)$ but has a 
different matrix realization.)

We will need the orthosymplectic supermatrices in the alternating 
format; in this format we take the matrix $B_{m, 2n}(\alt)=\antidiag 
(1, \dots, 1, -1, \dots, -1)$ with the only nonzero entries on the 
side diagonal, the last $n$ being $-1$'s.  The Lie superalgebra of 
such supermatrices will be denoted by $\fosp(\alt_{m|2n})$, where, as 
is easy to see, either $m=2n\pm 1$ or $m=2n$.

There is a 1-parameter family of deformations $\fosp_\alpha(4|2)$ of 
the Lie superalgebra $\fosp(4|2)$; its only explicit description we 
know (apart from \cite{BGLS}, of course) is in terms of Cartan matrix 
\cite{GL3}.

\ssec{4.1.  The superprincipal embeddings} Not every simple Lie 
superalgebra, even a finite dimensional one, hosts a superprincipal 
$\fosp(1|2)$-subsuperal\-geb\-ra.  Let us describe those that do.  
(Aside: an interesting problem is to describe {\it semiprincipal} 
embeddings into $\fg$, defined as the ones with the least possible 
number of irreducible components.)

We select the following basis in $\fosp(1|2)\subset \fsl(\bar 0|\bar 1|\bar 0)$: 
$$ 
X^-=\begin{pmatrix} 
0&0&0\\
0&0&0\\ 
-1&0&0\cr\end{pmatrix}, \; 
H=\begin{pmatrix} 
1&0&0\\0&0&0\\
0&  0&-1\\ \end{pmatrix}, \; 
X^+= \begin{pmatrix} 0& 0&1\\
0&0&0\\
0&0&0\\ 
\end{pmatrix}. 
$$
$$ 
\nabla^-=\begin{pmatrix} 0&0&0\\
1&0&0\\
0&1&0\\ 
\end{pmatrix}, \;
\nabla^+=\begin{pmatrix}
0&1&0\\
0&0&-1\\
0&0&0\\ \end{pmatrix}. 
$$
The highest weight $\fosp(1|2)$-module $\cM^\mu$ is illustrated with a graph
whose nodes correspond to the eigenvectors $l_i$ of $H$ with the weight indicated; the horisontal edges
depict the
$X^{\pm}$-action (the $X^+$-action is directed to the right, that of $X^-$ to
the left; each horizontal string is an irreducible $\fsl(2)$-submodule; two such
submodules are glued together into an $\fosp(1|2)$-module by the action of
$\nabla^{\pm}$ (we set $\nabla^+(l_n)=0$ and $\nabla^-(l_i)=l_{i-1}$; the corresponding
edges are not depicted below); we additionally assume that $p(l_\mu)=\bar 0$:   
$$ 
\begin{matrix} 
\dots\overset{\mu-2i}{\circ} \longleftrightarrow\overset{\mu-2i+2}{\circ}
\longleftrightarrow ...\quad ... \quad ...
\longleftrightarrow\overset{\mu-2}{\circ} \longleftrightarrow\overset{\mu}{\circ}\\
\dots\overset{\mu-2i+1}{\circ} \longleftrightarrow\overset{\mu-2i+3}{\circ}
\longleftrightarrow ...\longleftrightarrow\overset{\mu-3}{\circ}
\longleftrightarrow\overset{\mu-1}{\circ}\end{matrix}
$$ 
As follows from the relations of type 0 below in sec 4.2, the module $\cM^n$ for $n\in \Zee
_+$ has an irreducible submodule isomorphic to $\Pi(\cM^{-n-1})$; the quotient,
obvoiusly irreducible as follows from the same formulas,  will be denoted by $\cL^n$. 

Serganova completely described superprincipal embeddings of 
$\fosp(1|2)$ into a simple finite dimensional Lie superalgebra 
\cite{LSS} (the main part of her result was independently obtained in 
\cite{vJ}).

As the $\fosp(1|2)$-module corresponding to the superprincipal 
embedding, a simple finite dimensional Lie superalgebra $\fg$ is as 
follows (the missing simple algebras $\fg$ do not contain a 
superprincipal $\fosp(1|2)$):

\ssec{Table 4.1. $\fg$ that admits a superprincipal subalgebra: as the
$\fosp(1|2)$-module} 
$$
\renewcommand{\arraystretch}{1.4}
\begin{tabular}{|l|ll|l|}
\hline
$\fg$&$\fg=\cL^2\oplus
(\mathop{\oplus}\limits_{i>1} \cL^{2k_i})$ for $i\geq 2$&$
\oplus$&$(\mathop{\oplus}\limits_j\Pi(\cL^{m_j}))$ for $j\geq 1$\cr 
\hline
$\fsl(n|n+1)$&$\cL^2\oplus \cL^4 \oplus \cL^6 \dots \oplus \cL^{2n-2}$&&$\oplus
\Pi(\cL^1)\oplus \Pi(\cL^3)\oplus\dots \oplus \Pi(\cL^{2n-1})$\cr  
$\fosp (2n-1|2n)$&$\cL^2\oplus \cL^6 \oplus \cL^{10} \dots \oplus
\cL^{4n-6}$&&$\oplus  \Pi(\cL^3)\oplus \Pi(\cL^7)\oplus\dots \oplus
\Pi(\cL^{4n-1})$\cr    
($n>1$)&&&\cr 
$\fosp (2n+1|2n)$&$\cL^2\oplus \cL^6 \oplus \cL^{10}
\dots \oplus \cL^{4n-2}$&&$\oplus \Pi(\cL^3)\oplus \Pi(\cL^7)\oplus\dots
\oplus \Pi(\cL^{4n-1})$\cr  
$\fosp (2|2)\cong \fsl(1|2)$&$\cL^2$&&$\oplus \Pi(\cL^1)$\cr
$\fosp (4|4)$&$\cL^2\oplus \cL^6$&&$\oplus \Pi(\cL^3)\oplus \Pi(\cL^3)$\cr
$\fosp (2n|2n)$&$\cL^2\oplus \cL^6 \oplus \cL^{10} \dots \oplus 
\cL^{4n-2}$&$\oplus \cL^{2n-2}$&$\oplus \Pi(\cL^3)\oplus \Pi(\cL^7)\oplus\dots \oplus
\Pi(\cL^{4n-1})$\cr  
$\fosp (2n+2|2n)$&$\cL^2\oplus \cL^6 \oplus \cL^{10} \dots \oplus \cL^{4n+2}$&$\oplus
\cL^{2n}$&$\oplus \Pi(\cL^3)\oplus \Pi(\cL^7)\oplus\dots \oplus \Pi(\cL^{4n-1})$\cr  
$\fosp_{\alpha}(4|2)$&$\cL^2$&$\oplus \cL^2$&$\oplus  \Pi(\cL^3)$  \cr
\hline
\end{tabular} 
$$

The Lie superalgebra $\fg$ of type $\fosp$ that contains a 
superprincipal subalgebra $\fosp (1|2)$ can be generated by two 
elements.  For such elements we can take $X:=\nabla^+\in 
\cL^2=\fosp(1|2)$ and a lowest weight vector $Z:=l_{-r}$ from the 
module $M=\cL^r$ or $\Pi(\cL^r)$, where for $M$ we take $\Pi(\cL^3)$ 
if $\fg\not=\fosp(2n|2m)$ or the last module with the even highest 
weight vector in the above table (i.e., $\cL^{2n-2}$ if $\fg 
=\fosp(2n|2n)$ and $\cL^{2n}$ if $\fg =\fosp(2n+2|2n)$).

To generate $\fsl(n|n+1)$ we have to add to the above $X$ and $Z$ a 
lowest weight vector $U$ from $\Pi(\cL^1)$.  (Clearly, $Z$ and $U$ are 
defined up to factors that we can select at our convenience; we will 
assume that a basis of $L^r$ is fixed and denote $Z=t\cdot l_{-r}$ and 
$U=s\cdot l_{-1}$ for $t, s\in \Cee $.)

We call the above $X$ and $Z$, together with $U$, and fortified by 
$Y:=X^-\in L^2$ the {\it Jacobson's generators}.  The presence of $Y$ 
considerably simplifies the form of the relations, though slightly 
increases the number of them.

\ssec{4.2.  Relations between Jacobson's generators} We repeat the 
arguments from sec.  1.2.  Since we obtain the relations recurrently, 
it could happen that a relation of higher degree implies a relation of 
a lower degree.  This did not happen when we studied $\fsl(\lambda)$, 
but does happen in what follows, namely, relation 1.2 implies 1.1.

We divide the relations between Jacobson's generators into the types 
corresponding the number of occurence of $z$ in them: {\bf 0}.  
Relations in $\fsl(1|2)$ or $\fosp(1|2)$; {\bf 1}.  Relations coming 
from the $\fosp(1|2)$-action on $\cL^{2k_2}$; {\bf 2}.  Relations 
coming from $\cL^{2k_{1}}\wedge \cL^{2k_2}$; {\bf 3}.  Relations 
coming from $\cL^{2k_2}\wedge \cL^{2k_2}$; $\pmb\infty$.  Relation 
that shear the dimension.

The relations of type 0 are the well-known relations in $\fsl(1|2)$, 
those of them that do not involve $U$ (marked with an $*$) are the 
relations for $\fosp(1|2)$.  The relations of type 1 that do not 
involve $U$ express that the space $\cL^{2k_2}$ is the 
$\fosp(1|2)$-module with highest weight $2k_2$.  To simplify notations 
we denote: $Z_i=\ad X^iZ$ and $Y_i=\ad
X^iY$.   $$
\begin{matrix}
\pmb{0.1}^*.& [Y, Y_1]=0, &\pmb{0.2}^*.& [Y_2,  Y]=2Y,&\pmb{0.3}^*.& [Y_2, X]=-X,\\  
\pmb{0.4}.& [Y, U]=0,& \pmb{0.5}.& [U, U] = -2Y;&\pmb{0.6}.& [U, Y_1]=0,\\           
\pmb{0.7}.& [[X, X], [X,  U]]=0,& \pmb{0.8}.&  [Y_2, U]=U.\\ 
\end{matrix}
$$ 
$$
\pmb{1.1}.\; \; [Y, \, Z] = 0\Longleftarrow \pmb{1.2}.\; \; [[X, Y], \, Z] = 0, \quad
\pmb{1.3}.\; \; Z_{4k_{1}} = 0, \quad \pmb{1.4}.\; \; [Y_2, Z]= 3Z.
$$

\ssbegin{4.3}{Theorem}
For the Lie superalgebras indicated, all the relations between Jacobson's
generators are the above relations of types $0, 1$  and the relations from 
$\S 6$. \end{Theorem}

\section*{\S 5. The Lie superalgebra $\fgl (\lambda|\lambda+1)$ as the quotient of
$\fdiff (1|1)$ and a subalgebra of
$\fsl_+ (\infty|\infty)$}

There are several ways to superize $\fsl_+ (\infty|\infty)$.  For a 
description of \lq\lq the best" one from a certain point of view see 
\cite{E}.  For our purposes any version of $\fsl_+ (\infty|\infty)$ 
will do.

\ssec{5.1} The Poincar\' e-Birkhoff-Witt theorem states that 
$U(\fosp(1|2))\cong \Cee [X^-, \nabla^-, H, \nabla^+, X^+]$, as 
superspaces.  Set 
$U_\lambda=U(\fosp(1|2))/(\Delta-\lambda^2+\frac{9}{4})$.  
Denote: $\partial_x=\frac{\partial}{\partial x}$, $\partial_\theta= 
\frac{\partial}{\partial\theta}$ and set
$$
X^-=-\partial_x, \quad \nabla^-=\partial_\theta -\theta \partial_x, 
\quad H=2x\partial _x+\theta\partial_\theta (\lambda-1),\quad 
\nabla^+=x\partial_\theta +x\theta\partial _x-\lambda\theta, \quad 
X^+=x^2\partial_x-(\lambda-1) x.
$$ 
These formulas establish a morphism of $\fosp(1|2)$-modules and, 
moreover, of associative superalgebras: $U_\lambda\longrightarrow \Cee 
[x, \theta , \partial_x, \partial_\theta ]$.

In what follows we will need a well-known fact: the Casimir operator
$$
\Delta=2(X^+X^-+X^-X^+)+\nabla^+\nabla^--\nabla^-\nabla^++H^2
$$ 
acts on the irreducible $\fosp(1|2)$-module $\cL^{\mu}$ as the scalar 
operator of multiplication by $\mu^2+3\mu$.  (The passage from $\mu$ 
to $\lambda$ is done with the help of a shift by $\frac{3}{2}$.)

Consider the Lie superalgebra $LU(\fosp(1|2))$ associated with the 
associative superalgebra $U_\lambda$.  It is easy to see that, as 
$\fosp(1|2)$-module,
$$ 
LU_\lambda=\cL^0\oplus \cL^2\oplus\dots\oplus
\cL^{2n}\oplus \dots \oplus\Pi (\cL^1\oplus \cL^3\oplus \dots)\eqno{(5.1)}
$$

In the same way as for Lie algebras we show that $LU_n$ contains an ideal
$I_{n}$ for $n\in \Nee \setminus\{0\}$ and the quotient $LU_n/I_{n}$ is the
conventional $\fsl(n|n+1)$. It is clear that for $\lambda \neq \Zee $ the Lie
algebra $LU_\lambda$ has only one ideal --- the space $\cL^0$ of constants and
$LU_\lambda=\cL^0\oplus [LU_\lambda, LU_\lambda]$; hence, there is a supertrace
on $LU_\lambda$. This justifies the following notations    
$$ 
\fsl(\lambda|\lambda+1)=\fgl(\lambda|\lambda+1)/\cL^0,\quad
\text{ where }\fgl(\lambda|\lambda+1)=\left\{\begin{matrix} U_\lambda&\text{for
$\lambda \neq
\Nee\setminus\{0\}$}\\ 
LU_n/I_n&\text{otherwise.}\end{matrix}\right.\eqno{(5.2)}
$$
The definition directly implies that
$\fgl(-\lambda|-\lambda+1)\cong\fgl(\lambda|\lambda+1)$, so speaking about real
values of $\lambda$ we can confine ourselves to the nonnegative values.

Define $\fosp(\lambda+1|\lambda)$ as the Lie subsuperalgebra of
$\fsl(\lambda|\lambda+1)$ invariant with respect to the
involution 
$$
X\to \left\{\begin{matrix}
 -X&\text{if $X\in \cL^{4k}$ or $X\in \Pi(\cL^{4k\pm 1})$}\cr
X&\text{if $X\in \cL^{4k\pm 2}$ or $X\in \Pi(\cL^{4k\pm
3})$},\end{matrix}\right.\eqno{(5.3)}
$$
which is the analogue of the map 
$$
X\to -X^{st}\quad \text{for}\;  \quad X\in\fgl(m|n).\eqno{(5.4)}
$$

\ssec{5.2.  The Lie superalgebras $\fsl(*|*+1)$ and $\fosp(2*|*+1)$, 
for $*\in\Cee P^1=\Cee \cup\{*\}$} The \lq\lq dequantization" of the 
relations for $\fsl(\lambda|\lambda+1)$ and $\fosp(\lambda+1|\lambda)$ 
is performed by passage to the limit as $\lambda\longrightarrow\infty$ 
under the change $ t\mapsto\frac{t}{\lambda}$.  We denote the limit 
algebras by $\fsl(*|*+1)$ and $\fosp(*+1|*)$ in order not to mix them 
with $\fsl(\infty|\infty+1)$ and $\fosp(\infty|\infty+1)$, 
respectively.

\section*{\S 6. Tables. The Jacobson generators and relations between them}

\ssec{Table 6.1. Infinite dimensional case} 

$\bullet\quad \underline{\fosp (\lambda|\lambda+1)}$. Generators: 
$$
X = x \partial_\theta + x \theta \partial_x - \lambda \theta, \; \; 
Y = \partial_x,\; \;  
Z = t(\partial_x\partial_\theta - \theta {\partial_x}^2).
$$
Relations:
$$
\begin{tabular}{rl}
{\pmb 2.1}.& $3 [Z, Z_3] + 2 [Z_1, Z_2] = 6t(2\lambda+1) Z$, \cr {\pmb 
2.2}.& $[Z_1, Z_3] = 2 t^2(\lambda-1) (\lambda+2) Y +2t(2\lambda+1) 
Z_1$, \cr {\pmb 3.1}.& $[Z_1, [Z, Z]] = 0$.
\end{tabular}
$$

\noindent $\bullet\quad \underline{\fosp (*|*+1)}$.  Relations: the 
same as in sec 4.2 plus the following relations:
$$
\begin{tabular}{rl}
\pmb{2.1}.& $3[Z, Z_3]+2[Z_1, Z_2]=12t Z$,\cr 
\pmb{2.2}.&$[Z_1, Z_3]=2t^2Y+4tZ_1$.   
\end{tabular}
$$

\noindent $\bullet\quad \underline{\fsl (\lambda +1|\lambda)}$ for 
$\lambda\in \Cee P^1$.  Generators (for $\lambda\in \Cee $): the same 
as for $\fosp (\lambda|\lambda+1)$ and $U=\partial_\theta -\theta 
\partial_x$.

Relations: the same as for $\fosp (\lambda|\lambda+1)$ plus the 
following
$$
\begin{tabular}{c l c l}
\pmb{1.5}.\; &$3 [Z, [X, U]] - [U, Z_1] = 0$,&\pmb{2.3}.\; &$[Z, [U, Z]] = 0$,\cr
\pmb{1.6}.\; &$[[X, U], Z_1] = 0$,&\pmb{2.4}.\; &$[Z_1, [U, Z]] = 0$.\cr
\end{tabular}
$$

\ssec{Table 6.2. Finite dimensional algebras} 
In this table $E_{ij}$ are the matrix units; $X^{\pm}_i$ stand for the
conventional Chevalley generators of $\fg$. 
 
$\bullet\quad \underline{\fsl(n+1|n)}$ for $n\geq 3$. Generators:
$$
\begin{matrix}
X = \sum\limits_{1\leq i \leq n} 
	\bigl((n-i+1)E_{2i-1, 2i} - i E_{2i, 2i+1}\bigr),&  Y = \sum\limits_{1\leq i \leq
2n-1} E_{i+2, i},\\
U = \sum\limits_{1\leq i \leq 2n} (-1)^{i+1} E_{i+1, i},&  Z =
\sum\limits_{1\leq i \leq 2n-2}(-1)^{i+1} E_{i+3, i}.
\end{matrix}
$$
{\bf Relations}: those for $\fsl(\lambda+1|\lambda)$ with $\lambda=n$ and an extra
relation to shear the dimension: 
$$
(\ad Z)^n([X, X])=0.
$$

\noindent For $n=1$ the relations degenerate in relations of type 0.

$\bullet\quad \underline{\fosp (2n+1|2n)}$. Generators: 
$$
\begin{tabular}{rl}
$X$=&$ \sum\limits_{1 \leq i \leq n}\bigl((2n-i+1)(E_{2i-1, 
2i}+E_{4n+2-2i, 4n+3-2i}) -i(E_{2i, 2i+1}-E_{4n+1-2i, 
4n+2-2i})\bigr)$, \cr $Y$=&$ E_{2n+2, 2n} + \sum\limits_{1 \leq i \leq 
2n-1} (E_{i+2, i}-E_{4n+2-i, 4n-i})$, \cr $Z$&$- E_{2n+2, 2n-1} - 
E_{2n+3, 2n} +\sum\limits_{1 \leq i \leq 2n-2}\bigr((-1)^i E_{i+3, i} 
+ E_{4n+2-i, 4n-1-i}\bigr)$.
\end{tabular}
$$

{\bf Relations}: those for $\fosp(2\lambda+1|2\lambda)$ with 
$\lambda=n$ and an extra relation to shear the dimension (the form of 
the relation is identical to that for $\fsl(n+1|n)$).

$\bullet\quad \underline{\fosp_\alpha (4|2)}$.  Generators: As 
$\fosp(1|2)$-module, the algebra $\fosp_\alpha (4|2)$ has 2 isomorphic 
submodules.  The generators $X$ and $Y$ belong to one of them.  It so 
happens that we can select $Z$ from either of the remaining submodules 
and still generate the whole Lie superalgebra.  The choice (a) is from 
$\Pi(\cL^3)$; it is unique (up to a factor).  The choices (b) and (c) 
are from $\cL^2$; one of them seem to give simpler relations.
$$
\begin{tabular}{rl}
$X$&$-\frac{\alpha+1}{\alpha}X^+_1+\frac{\alpha}{\alpha+1}X^+_2+
\frac{1}{\alpha(\alpha+1)}X^+_3,\quad Y = [X^-_1, X^-_2]+[X^-_1, X^-_3]+[X^-_2,
X^-_3]$, \cr  
$Z$&$\left \{\begin{matrix} \text{a)}&-[[X^-_1, X^-_2], X^-_3]];\\
\text{b)}&-(1+2\alpha)[X^-_1, X^-_2]+\alpha^2(2+\alpha)[X^-_1, X^-_3]+
(\alpha -1)(1+\alpha)^2[X^-_2, X^-_3];\\
\text{c)}& -[X^-_2, X^-_3]-(\alpha+1)[X^-_1, X^-_2].\end{matrix}\right .$\cr
\end{tabular}
$$
Relations of type 0 are common for cases a) -- c):
$$
\pmb{0.1}\; [Y, [Y, [X, X]]]=4Y;\quad \quad \pmb{0.2}\; [Y_1 [X,
X]]]=-2X;
$$
The other relations are as follows.

{\bf Relations a)}:
$$
\begin{tabular}{rl}
\pmb{1.1}& $[Y_1, Z_1]=3Z, \quad\quad\quad\quad\quad \pmb{1.2}\; (\ad [X,
X])^3Z_1=0$;\cr
\pmb{2.1}& $[Z, Z]=0;\quad\quad\quad\quad\quad\quad\quad \pmb{2.2}\; [Z_1, [[X,
X], Z]]=-4\frac{\alpha^2+\alpha+1}{\alpha(\alpha+1)}Z$, \cr
\pmb{3.1}& $[\ad [X, X](Z_1), [Z_1, \ad [X, X](Z_1)]]= $\cr
&$-\frac{16}{\alpha(\alpha+1)}Y+8\frac{\alpha^2+\alpha+1}{\alpha(\alpha+1)}[Z_1, 
\ad [X, X](Z_1)]+16\frac{(\alpha^2+\alpha+1)^2}{\alpha^2(\alpha+1)^2}Z_1$.\cr
\end{tabular}
$$

{\bf Relations b)}:
$$
\begin{tabular}{rl}
\pmb{1.1}& $[Y_1, Z_1]=2Z; \quad \quad\quad\quad\quad \pmb{1.2}\; (\ad [X,
X])^2Z_1=0$;\cr
\pmb{2.1}*& $[Z_1, Z_1]=2[Z, [Z,[X, X]]]-18\alpha^2(1+\alpha)^2Y+
4(1-\alpha)(2+\alpha)(1+2\alpha)Z$;   \cr
\pmb{3.1}& $(\ad Z)^3X=0$,\cr
\pmb{3.2}*& $[[Z, Z_1], (\ad [X, X])^2Z_1]=$\cr
&$(-1+\alpha)(2+\alpha)(1+2\alpha)[Z,
[Z,[X, X]]] +12(1-\alpha)(2+\alpha)(1+2\alpha)\alpha^2(1+\alpha)^2Y+$\cr
&$8
(1-3\alpha^2-\alpha^3)(-1-3\alpha+\alpha^3)Z$.\cr
\end{tabular}
$$

{\bf Relations c)}: same as for b) except that the relations marked in b) by an ${}*$
should be replaced with the following ones
$$
\begin{tabular}{rl}
\pmb{2.1}& $[Z_1, Z_1]=2[Z, [Z,[X, X]]]-2\alpha^2Y+4(2+\alpha)Z$;  \cr 
\pmb{2.2}& $(\ad [X, X])^2Z_1=(-2-\alpha)[Z, [Z,[X, X]]] 
-8(1+\alpha)Z+4\alpha^2(2+\alpha)Y$.\cr
\end{tabular}
$$

\section*{\S 7. Remarks and problems}

\ssec{7.1.  On proof} For the exceptional Lie algebras and 
superalgebras $\fosp_\alpha(4|2)$ the proof is direct: the quotient of 
the free Lie algebra generated by $x, y$ and $z$ modulo our relations 
is the needed finite dimensional one.  For rank $\fg\leq 12$ we 
similarly computed relations for $\fg=\fsl(n)$, $\fo(2n+1)$ and 
$\fsp(2n)$; as Post pointed out, together with the result of \cite{PH} 
on deformation (cf.  2.7) this completes the proof for Lie algebras.  
The results of \cite{PH} on deformations can be directly extended for 
the case of $\fsl(2)$ replaced with $\fosp(1|2)$; this proves Theorem 
4.3.

Our Theorem 2.6 elucidates Proposition 2 of \cite{F}; we just wrote 
relations explicitely.  Feigin claimed \cite{F} that for $\fsl(\lambda)$ 
the relations of type 3 follow from the decomposition of 
$L^{2k_2}\wedge L^{2k_2}\subset L^{2k_2}\wedge L^{2k_2}\wedge 
L^{2k_2}$.  We verified that this is so not only in Feigin's case but 
for all the above-considered algebras except $\fe_6$, $\fe_7$ and 
$\fe_8$: for the latter one should consider the whole $L^{2k_2}\wedge 
L^{2k_2}\wedge L^{2k_2}$, cf.  \cite{GL1}.  Theorem 4.3 is a direct 
superization of Theorem 2.6.

\ssec{7.2.  Problems} 1) How to present $\fo(2n)$ and $\fosp(2m|2n)$?  
One can select $z$ as suggested in sec.  1.1.  Clearly, the form of 
$z$ (hence, relations of type 1) and the number of relations of type 3 
depend on $n$ in contradistinction with the algebras considered above.  
Besides, the relations are not as neat as for the above algebras.  We 
should, perhaps, have taken the generators as for $\fo(2n-1)$ and add 
a generator from $L^{2n-2}$.  We have no guiding idea; to try at 
random is frustrating, cf.  the relations we got for 
$\fosp_\alpha(4|2)$.

2) We could have similarly realized the Lie algebra $\fsl(\lambda)$ as 
the quotient of $U(\fvect(1))$, where $\fvect(1)=\fder\Cee [u]$.  
However, $U(\fvect(1))$ has no center except the constants.  What are 
the generators of the ideal --- the analog of (2.0) --- modulo which 
we should factorize $U(\fvect(1))$ in order to get $\fsl(\lambda)$?  
(Observe that in case $U(\fg)$, where $\fg$ is a simple finite 
dimensional Lie superalgebra such that $Z(U(\fg))$ is not noetherian, 
the ideal --- the analog of (5.0) --- is, nevertheless, finitely 
generated, cf.  \cite{GL2}.)

3) Feigin realized $\fsl(*)$ on the space of functions on the open 
cell of $\Cee P^1$, a hyperboloid, see \cite{F}.  Examples of 
\cite{DGS} are similarly realized.  Give any realization of 
$\fo/\fsp(*)$ and its superanalogs.

4) Other 
problems are listed in sec. 8.1--8.3 below.

\ssec{7.3. Serre relations are more convenient than our
ones}  The following Table represents results of V.~Kornyak's
computations. $N_{GB}$ is the number of relations in Groebner basis,
$N_{comm}$ is the number of non-zero commutators in the multiplication
table, $D_{GB}$ is a maximum degree of relations in $GB$, Space is measured in in bytes. The
corresponding values for Chevalley generators/Serre relations are given in
brackets. 

$$
\begin{matrix} 
\text{alg}& N_{GB}    &     N_{comm}  &D_{GB}& \text{Space}    
& \text{Time}\\
\fsl(3)  &  23 \; (24)  &  21 \; (21) &  9\; (4)  & 1300 \; (1188)  & 
<1 sec  \;   (<1 sec)\\
\fsl(4)  &  69 \; (84)   & 70 \; (60)  &17\; (6)   &3888 \; (3612)   &
<1 sec  \;   (<1 sec)\\
\fsl(5)  &193\; (218)  & 220\; (126)  &25\; (8)  &13556 \; (8716)   &
<1 sec  \;   (<1 sec)\\
\fsl(6)  & 444\; (473)  & 476\; (225)  &33 (10)  &34692\; (18088)    &
2 sec  \;   (<1 sec)\\
\fsl(7)  & 893\; (908)  & 937\; (363)  &41 (12)  &80272\; (33700)   &1
0 sec   \;   (1 sec)\\
\fsl(8)  &1615  (1594)  &1632\; (546)  &49 (14) &162128\; (57908)   &34
sec   \;   (3 sec)\\
\fsl(9)  &2705  (2614)  &2714\; (780)  &57 (16) &314056\; (93452)  &109
sec\; (6 sec)\\
\fsl(10)  &4263\; (4063)  &4138\; (1071)  &65\; (18) &534684\; (143456)  &336
sec (10 sec)\\
\fsl(11) &6405\; (6048) & 6224\; (1425)  &73\; (20) &921972\; (211428) &1058 sec
\; (19 sec)
\end{matrix} 
$$
For the other Lie algebras, especially exceptional ones, the 
comparison is even more unfavourable.  Nevertheless, for 
$\fsl(\lambda)$ with noninteger $\lambda$ there are only the Jacobson 
generators and we have to use them.

\section*{\S 8. Lie algebras of higher ranks. The analogs of
the exponents and $W$-algebras}

The following Tables 8.1 and 8.2 introduce the generators for the Lie 
algebras $U_\fg(\lambda)$ and the analogues of the exponents that 
index the generalized $W$-algebras (for their definition in the 
simplest cases from different points of view see \cite{FFr} and 
\cite{KM}; we will follow the lines of \cite{KM}).

Recall that (see 0.1) one of the definitions of $U_\fg(\lambda)$ is as 
the associative algebra generated by $\tilde \fg$; we loosely denote 
it by $\tilde S^{\bcdot}(\tilde \fg)$.  For the generators of 
$LU_\fg(\lambda)$ we take the Chevalley generators of $\fg$ (since by 
7.3 they are more convenient) and the lowest weight vectors of the 
irreducible $\fg$-modules that constitute $\tilde S^2(\tilde \fg)$.

\ssec{8.1.  The exponents} This section is just part of Table 1 from 
\cite{OV} reproduced here for the convenience of the reader.  Recall that 
if $\fg$ is a simple (finite dimensional) Lie algebra, $W=W_\fg$ is 
its Weyl group, $l=\rk~\fg$, $\alpha_1$, \dots , $\alpha_l$ the simple 
roots, $\alpha_0$ the lowest root; the $n_i$ the coefficients of 
linear relation among the $\alpha_i$ normed so that $n_0=1$; let 
$c=r_1\cdot \dots\cdot r_l$, where $r_i$ are the reflections from $W$ 
associated with the simple roots, be the Killing--Coxeter element.  
The order $h$ of $c$ (the Coxeter number) is equal to $\sum_{i>0} 
n_i$.  The eigenvalues of $c$ are $\varepsilon^{k_{1}}$, \dots, 
$\varepsilon^{k_{l}}$, where $\varepsilon$ is a primitive $h$-th root 
of unity.  The numbers $k_i$ are called the {\it exponents}.  Then

The exponents $k_i$ are the respective numbers $k_i$ from Table 1.1, 
e.g., $k_1=1$.  The number of roots of $\fg$ is equal to $l\sum_{i>0} 
n_i= 2\sum_{i>0} k_i$.  The order of $W$ is equal to
$$
zl!\prod n_i=\prod_{i>0} (k_i+1),
$$ 
where $z$ is the number of 1's among the $n_i$'s for $i>0$ (the number 
$z$ is also equal to the order of the centrum $Z(G)$ of the simply 
connected Lie group $G$ with the Lie algebra $\fg$).  The algebra of 
$W$-invariant polynomials on the maximal diagonalizable (Cartan) 
subalgebra of $\fg$ is freely generated by homogeneous polynomials of 
degrees $k_i+1$.
 
We will use the following notations: 

For a finite dimensional irreducible representations of finite 
dimensional simple Lie algebras $R(\lambda)$ denotes
the irreducible representation with highest weight $\lambda$ and
$V(\lambda)$ the space of this representation;

$\rho=\frac12\sum\limits_{\alpha>0}\alpha$ or $\rho$ is a weight such 
that $\rho(\alpha_i)=A_{ii}$ for each simple root $\alpha_i$.

The weights of the Lie algebras $\fo(2l)$ and $\fo(2l+1)$, $\fsp(2l)$ 
and $\ff_4$ ($l=4$) are expressed in terms of an orthonormal basis 
$\varepsilon_1$, \dots , $\varepsilon_l$ of the space $\fh^*$ over 
$\Qee$.  The weights of the Lie algebras $\fsl(l+1)$ as well as 
$\fe_7$, $\fe_8$ and $\fg_2$ ($l=7, 8$ and 2, respectively) are 
expressed in terms of vectors $\varepsilon_1$, \dots , 
$\varepsilon_{l+1}$ of the space $\fh^*$ over $\Qee$ such that $\sum 
\varepsilon_{i}=0$.  For these vectors $(\varepsilon_{i}, 
\varepsilon_{i})=\frac{l}{l+1}$ and $(\varepsilon_{i}, 
\varepsilon_{j})=\frac{1}{l+1}\quad \text{for } i\neq j$.  The indices 
in the expression of any weight are assumed to be different.

The analogues of the exponents for $LU_{\fg}(\lambda)$ are the 
highest weights of the representations that 
constitute $\tilde S^k(\tilde \fg)$.

\begin{Problem} Interprete these exponents in terms of the analog of 
the Weyl group of $LU_\fg(\lambda)$ in the sence of \cite{PS} and 
invariant polynomials on $LU_\fg(\lambda)$.
\end{Problem}

\ssec{8.2. Table. The Lie algebras $U_\fg(\lambda)$ as $\fg$-modules}

Columns 2 and 3 of this Table are derived from Table 5 in \cite{OV}.  
Columns 4 and 5 are results of a computer-aided study.  To fill in the 
gaps is a research problem, cf.  \cite{GL2} for the Lie algebras 
different from $\fsl$ type.

$$
\renewcommand{\arraystretch}{1.4}
\begin{tabular}{|c|c|c|c|c|}
\hline
$\fg$&$\ad$&$\tilde S^2(\tilde \fg)$&$\tilde S^3(\tilde \fg)$&$\tilde S^k(\tilde \fg)$\cr
\hline
$\fsl(2)$&$R(2\pi)$&$R(4\pi)$&$R(6\pi)$&$R(2k\pi)$\cr
\hline
$\fsl(3)$&$R(\pi_1+\pi_2)$&$R(2\pi_1+2\pi_2)$&$R(3\pi_1+3\pi_2)$&$R(k\pi_1+k\pi_2)$\cr
       &                
&$R(\pi_1+\pi_2)$&$R(2\pi_1+2\pi_2)$&$R((k-1)\pi_1+(k-1)\pi_2)$\cr
\hline
\end{tabular} 
$$
$$
\renewcommand{\arraystretch}{1.4}
\begin{tabular}{|c|c|c|c|c|}
\hline
$\fsl(4)$&$R(\pi_1+\pi_3)$&$R(2\pi_1+2\pi_3)$&$R(3\pi_1+3\pi_3)$ & \cr
         &                  &$R(\pi_1+\pi_3)$&$R(2\pi_1+2\pi_3)$  & \cr
         &                  &$R(2\pi_2)$         &$R(2\pi_1+\pi_2)$  & \cr
          &                  &                            &$R(\pi_2+2\pi_3)$&\cr
         &                  &                            &$R(\pi_1+\pi_3)$&\cr
         &                  &                            &$R(\pi_1+2\pi_2+\pi_3)$&\cr
\hline
$\fsl(n+1)$&$R(\pi_1+\pi_n)$&$R(2\pi_1+2\pi_n)$&$R(3\pi_1+3\pi_n)$
&\quad\quad \cr 
$n\geq 4$& &$R(\pi_1+\pi_n)$&$R(2\pi_1+2\pi_n)$& \cr 
& &$R(\pi_2+\pi_{n-1})$&$R(2\pi_1+\pi_{n-1})$& \cr 
& & &$R(\pi_2+\pi_{n-1})$&\cr 
& & &$R(\pi_2+2\pi_n)$&\cr 
& & &$R(\pi_1+\pi_n)$&\cr 
& & &$R(\pi_1+\pi_2+\pi_{n-1}+\pi_n)$&\cr 
\hline
\end{tabular} 
$$

The generators of $LU_{\fg}(\lambda)$ are the Chevalley generators 
$X^\pm_i$ of $\fg$ AND the lowest weight vectors from $\tilde S^2$.  
Denote the latter by $z_1$, $z_2$ (sometimes there is a third one, 
$z_3$).  Then the relations are (recall that $ h_i=[X^+_i, X^-_i]$):

(type 0) the Serre relations in $\fg$

(type 1) The relations between $X^\pm_i$ and $z_j$, namely:
$$
X^-_i(z_j)=0;\quad h_i(z_j)=\text{weight}_i(z_j);\quad (ad 
X^+_i)^{\text{the power determined by the weight of}z_j}(z_j)=0.
$$

\begin{Problem} Give an explicit form of the relations of higher
types. \end{Problem}

\ssec{8.3.  Tougher problems} Even if the explicit realization of the 
exceptional Lie algebras by differential operators on the base affine 
space were known at the moment, it is, nevertheless, a difficult 
computer problem to fill in the blank spaces in the above table and 
similar tables for Lie superalgebras.  To make plausible conjectures 
we have to compute $\tilde S^k(\tilde \fg)$ to, at least, $k=4$.  

Observe that for simple Lie algebras $\fg$ we have a remarkable 
theorem by Kostant which states that $U_{\fg}(\lambda)$ contains 
every finite dimensional irreducible $\fg$-module $V$ with 
multiplicity equal to the multiplicity of the zero weight in $V$; in 
view of which only the $\fsl(2)$-line is complete.

\section*{\S 9. A connection with integrable dynamical systems}

We recall the basic steps of the Khesin--Malikov construction and 
then superize them.

\ssec{9.1.  The Hamilton reduction} Let $(M^{2n}, \omega)$ be a 
symplectic manifold with an action $act$ of a Lie group $G$ on $M$ by 
symplectomorphisms (i.e., $G$ preserves $\omega$).  The derivative of 
the $G$-action gives rise to a Lie algebra homomorphism $a: 
\fg=Lie(G)\longrightarrow\fh(2n)$.  The action $act$, or rather, $a$ 
is called a {\it Poisson} one, if $a$ can be lifted to a Lie algebra 
homomorphism $\tilde a: \fg\longrightarrow\fpo(2n)$, where the Poisson 
algebra $\fpo(2n)$ is the nontrivial central extension of $\fh(2n)$.

For any Poisson $G$-action on $M$ there arises a $G$-equivariant map 
$p: M\longrightarrow\fg^*$, called the {\it moment map}, given by the 
formula
$$
\langle p(x), g\rangle=\tilde a(g)(x)\quad\text{for any}\quad x\in M, 
g\in\fg.
$$

Fix $b\in\fg^*$; let $G_b\subset G$ be the stabilizer of $b$.  Under 
certain regularity conditions (see \cite{Ar}) $p^{-1}(b)/G_b$ is a 
manifold.  This manifold is endowed with the symplectic form
$$
\begin{matrix}
\omega(\bar v, \bar w)=\omega(v, w)\quad\text{for arbitrary 
preimages}\quad v, w \;\text{of} \; \bar v, \bar w,\; 
\text{respectively}\\
\text{wrt the natural projection}\quad T(p^{-1}(b))\longrightarrow 
T(p^{-1}(b)/G_b).
\end{matrix} 
$$
The passage from $M$ to $p^{-1}(b)/G_b$ is called {\it Hamilton 
reduction}.  In the above picture $M$ can be the {\it Poisson 
manifold}, i.e., $\omega$ is allowed to be nondegenerate not on the 
whole $M$; the submanifolds on which $\omega$ is nondegenerate are 
called {\it symplectic leaves}.

{\it Example}.  Let $\fg=\fsl(n)$ and $M=\fg^*$, let $G$ be the group 
$N$ of uppertriangular matrices with 1 on the diagonal.  The coadjoint 
$N$-action on $\fg^*$ is a Poisson one, the moment map is the natural 
projection $\fg^*\longrightarrow\fn^*$ and $\fg^*/N$ is a Poisson 
manifold.

\ssec{9.2.  The Drinfeld--Sokolov reduction} Let $\fg=\hat \fa^{(1)}$, 
where $\fa$ is a simple finite dimensional Lie algebra (the case 
$\fa=\fsl(n)$ is the one considered by Gelfand and Dickii), hat 
denotes the Kac--Moody central extension.  The elements of $M=\fg^*$, 
can be identified with the $\fa$-valued differential operators:
$$
(f(t)dt, az^*)\mapsto (tf(t)+at\frac{d}{dt}))\frac{dt}{t}.
$$
Let $N$ be the loop group with values in the group generated by 
positive roots of $\fa$.  For the point $b$ above take the element 
$y\in\fa\subset\hat\fg^*$ described in \S 3.  If $\fa=\fsl(n)$, we can 
represent every element of $p^{-1}(b)/N$ in the form
$$
t\frac{d}{dt}+y+\begin{pmatrix} b_1(t)&\dots&b_n(t)\\
0&\dots&0\\
0&\dots&0\\
\end{pmatrix}\longleftrightarrow \frac{d^n}{d\varphi^n}+\tilde 
b_1(\varphi)\frac{d^{n-1}}{d\varphi^{n-1}}+\dots+\tilde b_n(\varphi),
$$

To generalize the above to $\fsl(\lambda)$, Khesin and Zakharevich 
described the Poisson--Lie structure on symbols of pseudodifferential 
operators, see \cite{KM} and refs therin.  Let us recall the main formulas.
 
\ssec{9.3.1. The Poisson bracket on symbols of $\Psi DO$} 
Set $D=\frac{d}{dx}$; define
$$
D^\lambda \circ f=fD^\lambda +\sum\limits_{k\geq 1}\binom{\lambda}{k} 
f^{(k)}D^{(\lambda-k)}, \quad\text{where}\; 
\binom{\lambda}{k}=\frac{\lambda(\lambda-1)\dots (\lambda-k+1)}{k!}.
$$
Set
$$
G_\lambda=\left\{D^\lambda (1+\sum\limits_{k\geq 1}u_k(x)D^{(-k)})\right\}
$$
and
$$
TG_\lambda=\left\{\sum\limits_{k\geq 1}v_k(x)D^{(-k)}\right\}\circ
D^\lambda,\quad T^* G_\lambda=D^{-\lambda}\circ DO.
$$
For $X=D^{-\lambda}\circ\sum\limits_{k\geq 0}u_k(x)D^{(k)}\in 
T^*G_\lambda$ and $L=\left(\sum\limits_{k\geq 
1}v_k(x)D^{(-k)}\right)\circ D^\lambda\in TG_\lambda $ define the 
pairing $\langle X, L\rangle$ to be
$$
\langle X, L\rangle=Tr(L\circ X),\quad\text{where}\quad Tr(\sum
w_k(x)D^{(k)})=\Res |_{x=0}w_{-1}.
$$
The Poisson bracket on the spae of psedodifferential symbols $\Psi 
DS_\lambda$ is defined on linear functionals by the formula
$$
\{X, Y\}(L)=X(H_Y(L)),\quad\text{where}\quad H_Y(L)=(LY)_+L-L(YL)_+.
$$

\begin{Theorem} {\em (Khesin--Malikov)} For $\fa=\fsl(\lambda)$ in the 
Drinfeld--Sokolov picture, the Poisson manifolds $p^{-1}(b)/N_b$ and 
$\Psi DS_\lambda$ are isomorphic.  Each element of the Poisson leaf 
has a representative in the form
$$
t\frac{d}{dt}+y+\begin{pmatrix} b_1(t)&\dots&b_n(t)&\dots\\
0&\dots&0&\dots\\
0&\dots&0&\dots\\
\end{pmatrix}\longleftrightarrow D^\lambda \left(1+\sum\limits_{k\geq 
1}\tilde b_k(\varphi)D^{(-k)}\right).
$$
\end{Theorem}

The Drinfeld--Sokolov construction \cite{DS}, as well as its 
generalization to $\fsl(\lambda)$ and $\fo{/}\fsp(\lambda)$ (\cite{KM}), 
hinges on a certain element that can be identified with the image of 
$X^+\in \fsl(2)$ under the principal embedding.  For the case of 
higher ranks this is the image in $U_\fg(\lambda)$ of the element 
$y\in \fg$ described in \S 3 for Lie algebras.  In $\fsl(\lambda)$ and 
$\fo{/}\fsp(\lambda)$ this image is just $\frac{d}{dx}$ (or the matrix 
whose only nonzero entries are the 1's under the main diagonal in the 
realization of $\fsl(\lambda)$ and $\fo{/}\fsp(\lambda)$ by matrices).

\ssec{9.4. Superization} 

\ssec{9.4.1 Basics} Further facts from Linear Algebra in Superspaces. 

The {\it tensor algebra} $T(V)$ of the superspace
$V$ is naturally defined:
$T(V)=\bigoplus_{n\geq 0} T^n(V)$, where $T^0(V)=k$ and
$T^n(V)=V\otimes \dots \otimes V$ ($n$ factors) for $n>0$.

The {\it symmetric algebra} of the superspace $V$ is $S(V)=T(V)/I$, 
where $I$ is the two-sided ideal generated by $v_1\otimes 
v_2-(-1)^{p(v_1)p(v_2)}v_2\otimes v_1$ for $v_1, v_2\in V$.

The {\it exterior algebra} of the superspace $V$ is $E(V)=S(\Pi(V))$.  
Clearly, both the exterior and symmetric algebras of the superspace 
$V$ are supercommutative superalgebras.  It is worthwhile to mention 
that if $V_{\ev }\neq 0$ and $V_{\od}\neq 0$, then both $E(V)$ and 
$S(V)$ are infinite dimensional.

A {\it Lie superalgebra} is defined with Sign Rule applied to the 
definition of a Lie algebra.  Its multiplication is called {\it 
bracket} and is usually denoted by $[\cdot, \cdot]$ or $\{\cdot , 
\cdot\}$.  If, however, we try to use this definition in attempts to 
apply the standard group-theoretical methods to differential equations 
on supermanifolds we will find ourselves at a loss: the supergroups 
and their modules are objects from different categories!  Accordingly, 
the following (equivalent to the conventional, ``naive'' one, see 
\cite{L1}) 
definition becomes useful: a {\it Lie superalgebra} is a superalgebra 
$\fg$ (defined over a field or, more generally, a supercommutative 
superalgebra $k$); the bracket should satisfy the following 
conditions: $[X, X]=0$ and $[Y, [Y, Y]]=0$ for any $X\in (C\otimes 
\fg)_{\ev }$ and $Y\in (C\otimes \fg)_{\od}$ and any supercommutative 
superalgebra $C$ (the bracket in $C\otimes \fg$ is defined via 
$C$-linearity and Sign Rule).
 
With an associative (super)algebra $A$ we associate Lie 
(super)algebras (1) $A_L$ with the same (super)space $A$ and the 
multiplication $(a, b)\mapsto [a, b]$ and (2) $\fder A$, the algebra 
of derivations of $A$, defined via Sign and Leibniz Rules.

From a Lie superalgebra $\fg$ we construct an associative superalgebra 
$U(\fg)$, called the {\it universal enveloping algebra} of the Lie 
superalgebra $\fg$ by setting $U(\fg)= T(\fg)/I$, where $I$ is the 
two-sided ideal generated by the elements $x\otimes 
y-(-1)^{p(x)p(y)}y\otimes x-[x, y]$ for $x, y\in \fg$.

The {\it Poincar\'{e}--Birkhoff--Witt theorem} for Lie algebras 
extends to Lie superalgebras with the same proof (beware Sign Rule) 
and reads as follows:

{\it if $\{X_i\}$ is a basis in $\fg_{\ev }$ and $\{Y_j\}$ is a basis 
in $\fg_{\od}$, then the monomials $X^{n_1}_{i_1}\dots 
X^{n_r}_{i_r}Y^{\varepsilon_1}_{j_1}\dots Y^{\varepsilon_s}_{j_s}$, 
where $n_i\in \Zee^+$ and $\varepsilon_j = 0, 1$, constitute a basis 
in the space $U(\fg)$.}

A superspace $M$ is called a {\it left module} over a superalgebra $A$ 
(or a {\it left $A$-module}) if there is given an even map {\it act}: 
$A\otimes M\rightarrow M$ such that if $A$ is an associative 
superalgebra with unit, then $(ab)m=a(bm)$ and $1m=m$ and if $A$ is a 
Lie superalgebra, then $[a, b]m=a(bm)-(-1)^{p(a)p(b)}b(am)$ for any 
$a, b \in A$ and $m\in M$.  The definition of a {\it right $A$-module} 
is similar.

\begin{rem*}{Convention} We endow every module $M$ over a 
supercommutative superalgebra $C$ with a two-sided module structure: 
the left module structure is recovered from the right module one and 
vice versa according to the formula $cm=(-1)^{p(m)p(c)}mc$ for any 
$m\in M$ and $c\in C$.  Such modules will be called $C$-{\it modules}.  
(Over $C$, there are {\it two} canonical ways to define a two-sided 
module structure, see \cite{L2}; the meaning of such an abundance is 
obscure.)
\end{rem*}

The functor $\Pi$ is, actually, tensoring by $\Pi(\Zee)$.  So there 
are two ways to apply $\Pi$ to $C$-modules: to get $\Pi 
(M)=\Pi(\Zee)\otimes_{\Zee}M$ and $(M)\Pi=M\otimes_{\Zee}\Pi(\Zee)$.  
The two-sided module structures on $\Pi (M)$ and $(M)\Pi$ are given 
via Sign Rule.

Sometimes, instead of the map {\it act} a morphism $\rho : A 
\rightarrow \End M$ is defined if $A$ is an associative superalgebra 
(or $\rho : A\longrightarrow (\End M)_L$ if $A$ is a Lie 
superalgebra); $\rho$ is called a {\it representation} of $A$ in $M$.

The simplest (in a sense) modules are those which are {\it 
irreducible}.  We distinguish {\it irreducible modules of $G$-type} 
(general); these do not contain invariant subspaces different from $0$ 
and the whole module; and their \lq\lq odd" counterparts, {\it 
irreducible modules of $Q$-type}, which do contain an invariant 
subspace which, however, is not a subsuperspace.  Consequently, {\it 
Schur's lemma} states that {\sl over $\Cee$ the centralizer of a set 
of irreducible operators is either $\Cee$ or} $\Cee \otimes \Cee 
^s=Q(1; \Cee)$, see the definition of the superalgebras $Q$ below.

The next in terms of complexity are {\it indecomposable} modules, 
which cannot be split into the direct sum of invariant submodules.

A $C$-module is called {\it free} if it is isomorphic to a module of 
the form $C\oplus \dots \oplus C\oplus\Pi (C)\oplus \dots \oplus\Pi 
(C)$ ($C$ occurs $r$ times, $\Pi (C)$ occurs $s$ times).

The {\it rank} of a free $C$-module $M$ is the element $\rk 
M=r+s\varepsilon$ from the ring $\Zee [\varepsilon]/ 
(\varepsilon^2-1)$.  Over a field, $C=k$, we usually write just $\dim 
M=(r, s)$ or $r|s$ and call this pair the {\it superdimension} of $M$.

The module $M^*=\Hom_C(M, C)$ is called {\it dual} to a $C$-module 
$M$.  If $(\cdot , \cdot)$ is the pairing of modules $M^*$ and $M$, 
then to each operator $F\in \Hom_C(M, N)$, where $M$ and $N$ are 
$C$-modules, there corresponds the dual operator $F^*\in \Hom_C(N^*, 
M^*)$ defined by the formula
$$
(F(m), n^*)=(-1)^{p(F)p(m)}(m, F^*(n^*)) \;\text{ for any }\; m\in M, 
\ n^*\in N^*.
$$
Over a supercommutative superalgebra $C$ a {\it supermatrix} is a 
supermatrix with entries from $C$, the parity of the matrix with only 
$(i, j)$-th nonzero element $c$ is equal to $p_(i)+p_(j)+p(c)$.  
Denote by $\Mat(\Par; C)$ the set of $\Par\times \Par$ matrices with 
entries from a supercommutative superalgebra $C$.

The even invertible elements from $\Mat(\Par; C)$ constitute the {\it 
general linear group} $GL(\Par ; C)$.  Put $GQ(\Par; C)=Q(\Par; C)\cap 
GL(\Par; C)$.

On the group $GL(\Par;C)$ an analogue of the determinant is defined; 
it is called the {\it Berezinian} (in honour of F.~A.~Berezin who 
discovered it).  In the standard format the explicit formula for the 
Berezinian is:
$$
\Ber\begin{pmatrix}
A\ B\\ 
D\ E\end{pmatrix}=\det (A-BE^{-1}D)\det E^{-1}.
$$
For the matrices from $GL(\Par; C)$ the identity $\Ber\ XY=\Ber\ 
X\cdot \Ber\ Y$ holds, i.e., $Ber: GL(\Par; C)\rightarrow GL(1|0; 
C)=GL(0|1; C)$ is a group homomorphism.  Set $SL(\Par; C)=\{X \in 
GL(\Par; C): Ber X=1\}$.  The {\it orthosymplectic} group of 
automorphisms of the bilinear form with the even canonical matrix is 
denoted (in the standard format) by $Osp(n|2m; C)$.

\ssec{9.4.2.  Pseudodiferential operators on the supercircle.  
Residues} Let $V$ be a superspace.  For $\theta=(\theta_1, \dots , 
\theta_n)$ set
$$
\begin{matrix}
V[x, \theta]=V\otimes \Kee [x, \theta];\; V[x^{-1}, x, \theta]=V\otimes
\Kee [x^{-1}, x, \theta];\\
V[[x^{-1}, \theta]]=V\otimes \Kee [[x^{-1},
\theta]];\\ 
V((x, \theta))=V\otimes \Kee [[x^{-1}]][x, \theta].\end{matrix}
$$
We call $V((x, \theta))x^\lambda$ the space of {\it pseudodifferential
symbols}. Usually, $V$ is a Lie (super)algebra. Such symbols correspond to
pseudodifferential operators (pdo) of the form 
$$
\sum\limits_{i=-\infty}^n\sum\limits_{k_{0}+\dots+k_{n}=i}
a_i(\partial_x)^{k_{0}}\theta_1^{k_{1}}\dots \theta_n^{k_{n}},
$$
Here $k_i=0$ or 1 for $i>0$ and $a_i(x, \theta)\in V$. This is clear. 

For any $P=\sum\limits_{i\leq m} P_ix ^i\theta_0^k\theta^j\in V((x, \theta))$
we call 
$P_{+}=\sum\limits_{i, j, k\geq 0} P_ix ^i\theta_0^k\theta^j$ the {\it
differential part} of $P$ and
$P_{-}=\sum\limits_{i, k< 0} P_ix ^i\theta_0^k\theta^j$ the {\it integral}
part of $P$. 

The space $\Psi DO$ of pdos is, clearly, the left module over the 
algebra $\cF $ of functions.  Define the left $\Psi DO$-action on $\cF 
$ from the Leibniz formula thus making $\Psi DO$ into a superalgebra.

Define the involution in the superalgebra $\Psi DO$ setting
$$
(a(t, \theta)D^i\tilde D^j)^*=(-1)^{jip(\tilde D)p(D)}\tilde D^jD^ia^*(x,
\theta).
$$

The following fact is somewhat unexpected.  If $D$ is an odd 
differential operator, then $D^2$ is well-defined as $\frac{1}{2}[D, 
D]$.  Hence, we can consider the set $V((x, \theta))x^\lambda$ for an 
{\it odd} $x$!  Therefore, there are two types of pdos: {\it contact} 
ones, when $D^2\neq 0$ for odd $D$'s and general ones, when all odd 
$D$'s are nilpotent.

For the definition of distinguished stringy superalgebgras crucial in 
what follows see \cite{GLS1}.

\begin{rem*}{Conjecture} There exists a residue for all distinguished 
dimensions, i.e. for the contact type pdos in dimensions $1|n$ for 
$n\leq 4$ and for the general pdos in dimensions $1|n$ for $n\leq 2$.
\end{rem*}

So far, however, the residue was defined only for contact type pdos, of
$\fk^L$ type, and only for $n=1$ at that. 

Let us extend \cite{MR} and define the {\it residue} of $P=\sum\limits_{i\leq m} 
P_ix ^i\theta_0^k\theta^j\in V((x, \theta_0, \theta))$ for $n=1$.  We 
can do it thanks to the following exceptional property of $\fk^L(1|1)$ 
and $\fk^M(1|1)$.  Indeed, over $\fk^L(1|1)$, the volume form \lq\lq 
$dt\pder{\theta}$" is, more or less, $d\theta$: consider the quotient 
$\Omega^1/\cF\alpha$, where $\alpha$ is the contact form preserved by 
$\fk^L(1|1)$; similarly, over $\fk^M(1|1)$, the transformation rules 
of $dt\pder{\theta}$" and $\tilde\alpha$, where $\tilde\alpha$ is the 
contact form preserved by $\fk^M(1|1)$, are identical.  
Therefore, define the residue by the formula
$$
\Res~P= \text{ coefficient of }\frac{\theta}{x}\text{
in the expansion of } P_{-1}.
$$

\begin{rem*}{Remark} Manin and Radul \cite{MR} considered the 
Kadomtsev--Petviashvili hierarchy associated with $\fns$, i.e., for 
$D=K_\theta$.  The formula for the residue allows one to directly 
generalize their result and construct a simialr hierarchy associated 
with $\fr$, i.e., for $D=\tilde K_\theta$.
\end{rem*}

This new phenomenon --- an invertible odd symbol --- doubles the old 
picture: let $\theta_0$ be the symbol of $D$, and let $x$ be the 
symbol of the differential operator $D^2$.  We see that the case of 
the odd $D$ reduces to either $V((x, \theta_0, \theta))x^\lambda$ or 
$V((x, \theta_0^{-1}, \theta))x^\lambda$.

\ssec{9.5.  Continuous Toda lattices} Khesin and Malikov \cite{KM} 
considered straightforward generalizations of the Toda lattices --- 
the dynamical systems on the orbits of the coadjoint representation of 
a simple finite diminsional Lie group $G$ defined as follows.  Let 
$\cX$ be the image of $X^+\in\fsl(2)$ in $\fg=Lie(G)$ under the 
principal embedding.  Having identified $\fg$ with $\fg^{*}$ with the 
help of the invariant nondegenerate form, consider the orbit 
$\cO_{\cX}$.  On $\cO_{\cX}$, the traces $H_{i}(A)=\tr(A+\cX)^i$ are 
the commuting Hamiltonians.

In our constructions we only have to consider in $LU_{\fg}(\lambda)$ 
either (for Lie algebras) the image of $\cX$ or (for Lie 
superalgebras) the image of $\nabla^+\in\fosp(1|2)$ under the 
superprincipal embedding of $\fosp(1|2)$.  For superalgebras we also 
have to replace trace with the supertrace.

For the general description of dynamical systems on the orbits of the 
coadjoint representations of Lie supergroups see \cite{LST}.  A 
possibility of odd mechanics is pointed out in \cite{LST} and in the 
subsequent paper by R.~Yu.~Kirillova in the same Procedings.  To take 
such a possibility into account, we have to consider analogs of the 
principal embeddings for $\fsq(2)$.  This is a full-time job; its 
results will be considered elsewhere.


\begin{thebibliography}{9999}
	
\bibitem[Ar]{Ar}
Arnold V., {\em Mathematical methods of classical mechanics}, 
Springer, 1989 
\bibitem[BO]{BO}
Yu.  A. Bakhturin, A. Yu.  Olshansky, The approximations and 
characteristic subalgebras of free Lie algebras, In: {\it Proc.  I. G. 
Petrovsky seminar}, 2, 1976, 145--150 (in Russian)
\bibitem[BWV]{BWV} 
Bergshoeff E., de Wit B., Vasiliev M., The structure of the 
super-$W_{\infty}(\lambda)$ algebra. Nucl. Phys B366, 1991, 315--346
\bibitem[BMP]{BMP}
Bouwknegt P., McCarthy J., Pilch K., Quantum group structure in the 
Fock space resolutions of $\widehat{\fsl(n)}$ representations, Commun.  
Math.  Phys., 131, 1990, 339--368
\bibitem[B]{B}
Burd\'ik \v C., Realizations of real semisimple Lie algebras: a method 
of construction, J. Phys.  A: Math.  Gen., 18, 1985, 3101--3111
\bibitem[BGLS]{BGLS}
Burd\'ik \v C., Grozman P., Leites D., Sergeev A., Realization of simple 
Lie algebras and superalgebras via creation and annihilation 
operators. I., Realization of Lie
algebras and superalgebras via creation and annihilation operators. 
I. Theor.  and Math.  Physics, v.  124, 2000, no.  2, 227--238;
English translation: 1048--1058

\bibitem[Di]{Di} 
Dixmier J. {\it Alg\` ebres envellopentes}, Gautier--Villars, Paris, 1974
\bibitem[DGS]{DGS}  
Donin J., Gurevich D., Shnider S., Quantization of function algebras 
on semisimple orbits in $\fg^{*}$.  Talk at Internat.  Conf.  
``Quantum groups and integrable systems'' June 19 -- 21, 1997, Prague.  
see also: Invariant quantization in one and two parameters on 
semisimple coadjoint orbits of simple Lie groups, math.QA/9807159
\bibitem[DS]{DS} 
Drinfeld V., Sokolov V., Lie algebras and equations of Korteveg--de 
Vries type, Current problems in mathematics, Itogi Nauki i Tekhniki, 
vol.~24, Akad.  Nauk SSSR, Vsesoyuz.  Inst.  Nauchn.  i Tekhn.  
Inform., 1984, pp.~81--180 (in Russian; for English translation see 
JOSMAR)
\bibitem[D]{D}
Dynkin E. B., Semi-simple subalgebras of semi-simple Lie algebras, 
Mat.  Sbornik, 30, 1952, 111--244 (AMS Transl., v.6, ser.  2, 1957)
\bibitem[E]{E}
Egorov G., How to superize $\fgl(\infty)$.  In: Mickelsson J., 
Peckonnen O. (eds.)  {\it Diff.  Geometric Methods in Theoretical 
Physics} (Proc.  conf 1991, Turku, Finland) World Sci., 1992, 135--146
\bibitem[F]{F}
Feigin B. L., The Lie algebras $\fgl (\lambda)$ and cohomologies of 
Lie algebra of differential operators, Russian Math.  Surveys, v.  43, 
2, 1988, 157--158
\bibitem[FFr]{FFr}
Feigin B. L., Frenkel E., Integrals of motion and quantum groups.  In: 
Donagi R. e.a.  (eds) {\it Integrable systems and quantum groups} LN 
in Math 1620, 1996, 349--418
\bibitem[FO]{FO}
Feigin B. L., Odessky A., Elliptic Sklyanin algebras, Funkt. Anal. 
Appl., 1989, v.23, n. 3, 45--54
\bibitem[Go]{Go} 
Golod P., A deformation of the affine Lie algebra $A_1^{(1)}$ and 
hamiltonian systems on the orbits of its subalgebras. In: {\it 
Group-theoretical methods in physics} Proc. of the 3rd seminar, 
Yurmala, 1985, v.1, Moscow, Nauka, 1986, 368--376
\bibitem[GL1]{GL1} 
Grozman P., Leites D., Defining relations associated with the 
principal $\fsl(2)$-subalgebras.  In: Dobrushin R., Minlos R., Shubin 
M., Vershik A. (eds.)  {\it Contemporary mathematical physics} (F. A. 
Berezin memorial volume), Amer.  Math.  Soc.  Transl.  Ser.2, vol.  
175, Amer.  Math.  Soc., Providence, RI, 1996, 57--68
\bibitem[GL2]{GL2} 
Grozman P., Leites D., Lie superalgebras of 
supermatrices of complex size: a closer view (to appear)

\bibitem[GL3]{GL3} 
Grozman P., Leites D., Defining relations for Lie superalgebras with 
Cartan matrix, hep-th 9702073; Czech.  J. Phys., Vol.  51,
2001, No.  1, 1--22

\bibitem[GLS1]{GLS1}
Grozman P., Leites D., Shchepochkina I., Lie superalgebras of 
string theories, hep-th 9702120

\bibitem[GLS2]{GLS2}
Grozman P., Leites D., Sergeev A., to appear

\bibitem[KS]{KS}
Kac V., Some algebras related to the quantum field theory, XI-th
All-Union Algebr.  Coll., Kishinev, 1971, 140--141 (in Russian);
Classification of simple Lie superalgebras.  Funkcional.  Anal.  i
prilozheniya 9, n.  3, 1975, 91--92; id., Letter to editors,
Funkcional.  Anal.  i prilozheniya 10, n.  2, 1976, 93; id., Lie
superalgebras, Adv.  Math., 1976, 55--110

\bibitem[K]{K}
Kac V., {\it Infinite dimensional Lie algebras}, 3rd edition,
Cambridge Univ.  Press, Cambridge, 1991

\bibitem[KR]{KR}
Kac V., Radul A., Quasifinite highest weight modules over the Lie 
algebra of differential operators on the circle, Commun.  Math.  
Phys., v.  157, 1993, 429--457

\bibitem[Ka]{Ka}
Kashivara M., Representation theory and $D$-modules on flag varieties, 
Asterisque, 173--174, 1989, 55--110
\bibitem[KM]{KM}
Khesin B., Malikov F., Universal Drinfeld--Sokolov reduction and the 
Lie algebras of matrices of complex size, Comm.  Math.  Phys., v.  
175, 1996, 113--134

\bibitem[KV]{KV} 
Konstein S., Vasiliev M., Supertraces on the algebras of observables 
of the rational Calogero model with harmonic potential. J. Math. Phys 
37(6), 1996, 2872--2891

\bibitem[L1]{L1}
Leites D., Quantization and supermanifolds.  In: F. Berezin, M. Shubin 
{\it Schr\"odinger equation}, Kluwer, Dordrieht, 1991
\bibitem[L2]{L2}
Leites D. (ed.)  {\it Seminar on supermanifolds}.  Reports of 
Stockholm Univ., $\#\#$1--34, 1987--92

\bibitem[LM]{LM}
Leites D., Montgomery S., New simple filtered Lie superalgebras: a 
construction (to appear)

\bibitem[LP]{LP}
Leites D., Poletaeva E., Defining relations for classical Lie algebras 
of polynomial vector fields.  Math.  Scand., 81, 1997, no. 1, 5--19

\bibitem[LSS]{LSS}
Leites D., Saveliev M., Serganova V., Embeddings of $\fosp(n|2)$ and 
the associated nonlinear supersymmetric equations.  In: (Markov, M.A., 
Man'ko, V.I. and Dodonov, V.V., eds.)  {\it Group Theoretical Methods 
in Physics} (Yurmala 1985), VNU Science Press, Utreht, v.1, 1986, 
255--297

\bibitem[LST]{LST}
Leites D. , Semenov-Tian-Shansky M., Integrable systems and Lie 
superalgebras.  In: L.~D.~Faddeev (ed.  ) Differential geometry, Lie 
groups and mechanics, V. Zap.  Nauchn.  Sem.  Leningrad.  Otdel.  Mat.  
Inst.  Steklov.  (LOMI), Nauka, Leningrad, v.  123, 1983, 92--97; 
Kirillova, R. Yu.  Explicit solutions of superized Toda lattices.  
ibid, 123 ,1983, 98--111

\bibitem[LS]{LS}
Leites D., Serganova, V., Defining relations for classical Lie 
superalgebras.  I., In: Mickelsson J., Peckonnen O. (eds.)  {\it Diff.  
Geometric Methods in Theoretical Physics} (Proc.  conf 1991, Turku, 
Finland) World Sci., 1992, 194--201

\bibitem[LAS]{LAS}
Leites D., Sergeev A., Orthogonal polynomials of discrete variable and
Lie algebras of complex size matrices (with A. Sergeev) In: Procedings
of M.~Saveliev memorial conference, MPI, Bonn, February, 1999,
MPI-1999-36, 49--70; Theor.  and Math.  Physics, v.  123, no.  2,
205--236 (Russian), 582--609


\bibitem[LSc]{LSc}
Leites D., Shchepochkina I., Classification of simple Lie superalgebras of
vector fields (to appear) 

\bibitem[LSc1]{LSc1} 
Leites D., Shchepochkina I., How to quantize
antibracket, Theor.  and Math.  Physics, to appear

\bibitem[MR]{MR}
Manin Yu., Radul A., A supersymmetric extension of the 
Kadomtsev--Petviashvili hierarchy.  Comm.  Math.  Phys.  98 (1985), 
no.  1, 65--77

\bibitem[M]{M}
Montgomery S., Constructing simple Lie superalgebras from associative 
graded algebras, J. Algebra 195 (1997), no. 2, 558--579

\bibitem[OV]{OV} 
Onishchik A. L., Vinberg \'E. B., {\it Seminar on algebraic groups and 
Lie groups}, Springer, 1990

\bibitem[PS]{PS} 
Penkov I., Serganova V., Generic irreducible representations of finite 
dimensional Lie superalgebras, International J. Math., v.  5, 1994, 
389--419
\bibitem[PH]{PH}  Post G., Hijligenberg N. van den, $\fgl(\lambda)$
and differential operators preserving polynomials. Acta Appl. 
Math., v. 44, 1996,  257--268

\bibitem[SV]{SV} 
Saveliev M., Vershik A., Continuum analogues of contragredient Lie 
algebras, Commun. Math. Phys., 126, 1989, 367--378

\bibitem[ShV]{ShV} 
Shoihet B., Vershik A., Graded lie algebras whose Cartan subalgebra is
the algebra of polynomials in one variable, Teoreticheskaya i
Matematicheskaya fizika, v.  123, 2000, no.  2, 345--352, English
translation in Theor.  and Mathem.  Phys., v.  123, 2000, no.  2,


\bibitem[SH]{SH} 
Sheinman O.M., Heighest weight modules over certain quasigraded Lie 
algebras over elliptic curves, Funct. Anal. Appl., 26: 3, 1992, 
203--208

\bibitem[S]{S} 
Sergeev A., Invariant polynomials on Lie superalgebras.  In [L2], v. 
32.; for details see also The invariant polynomials on simple Lie
superalgebras, math.RT/9810111; Represent.  Theory 3 (1999), 250--280;
Leites D., Sergeev A., Casimir operators for Poisson Lie
superalgebras.  In: Ivanov E. et.  al.  (eds.)  {\it Supersymmetries
and Quantum Symmetries} (SQS'99, 27-31 July, 1999), Dubna, JINR, 2000,
409--411

\bibitem[Sh]{Sh} 
Shoikhet B., Certain topics on the Lie algebra $\fgl(\lambda)$ 
representation theory, q-alg/9703029; Complex analysis and
representation theory, 1.  J. Math.  Sci.  (New York) 92 (1998), no. 
2, 3764--3806

\bibitem[vJ]{vJ} 
Van der Jeugt J., Principal five-dimensional subalgebras of Lie 
superalgebras, J. Math.  Phys, v.27 (12), 1986, 2842--2847

\end{thebibliography}
\end{document}